\documentclass[a4paper,11pt]{article}

\RequirePackage{amsmath,amssymb,amsthm,amsfonts}
\RequirePackage[numbers]{natbib}

\theoremstyle{plain}
\newtheorem{theorem}{Theorem}[section]
\newtheorem{proposition}[theorem]{Proposition}
\newtheorem{lemma}[theorem]{Lemma}
\newtheorem{corollary}[theorem]{Corollary}

\theoremstyle{definition}
\newtheorem{definition}[theorem]{Definition}

\newtheorem{remark}[theorem]{Remark}

\newtheorem{assumption}{Assumption}

\begin{document}

\begin{center}

  {\Large\bf
    Discrete Gaussian free fields on Hamming graphs}
  
  \normalsize

  \bigskip By \bigskip

  \textsc{Shuhei Mano}

  \smallskip
  
  The Institute of Statistical Mathematics, Japan
  
  %\bigskip
  %November 18, 2025-\\  
  %Version: \today

\end{center}

\small

{\bf Abstract.}
  We discuss a class of discrete Gaussian free fields
  on Hamming graphs
  of dimension $d$ with vertex set
  $\mathcal{X}=\{0,1,\ldots,n-1\}^d$,
  where interactions depend solely on the Hamming
  distance between vertices. This class provides
  a setting distinct from the commonly studied
  fields on the integer lattice with nearest-neighbour
  interactions.
  Using spectral and combinatorial properties of
  Hamming graphs,
  we show that the field $(g_x)_{x\in\mathcal{X}}$
  converges weakly to i.i.d.\ Gaussian random
  variables in the limit as $d\to\infty$ or
  $n\to\infty$, while the scaled field
  $(mg_x)_{x\in\mathcal{X}}$ becomes perfectly
  correlated in the massless limit $m\to 0$.
  We then investigate several specific models.
  For a model in which the interaction weights
  follow a binomial probability mass function
  as a function of Hamming distance, with parameter
  $\gamma\in(0,1-1/n]$, coupling $\gamma$ to $m$
  yields a double-scaling limit in which
  $(mg_x)_{x\in\mathcal{X}}$ has correlations
  that are convex in Hamming distance.
  Remarkably, at every distance, the correlation is
  exactly the reciprocal of the number of vertices
  at that distance from a given vertex, attaining
  its minimum at a distance strictly between $0$ and $d$.
  
\smallskip


\smallskip


\normalsize

\section{Introduction}
\label{sect:intr}

Consider a finite set $\mathcal{X}=\{0,1,\ldots,n-1\}^d$,
where $d,n\ge 2$. Set $\rho(x,x')=|\{i:x_i\neq x'_i\}|$
for $x=(x_i)\in\mathcal{X}$ and $x'=(x'_i)\in \mathcal{X}$,
which is called the {\it Hamming distance}. The Hamming
distance induces a relation on $\mathcal{X}\times \mathcal{X}$
by
\[
  (x,x')\in R_i ~ \Leftrightarrow ~ \rho(x,x')=i \quad
  \text{for} ~ i\in\{0,1,\ldots,d\}.
\]
An undirected graph $(\mathcal{X},R_1)$ with vertices
$\mathcal{X}$ and edges $R_1$ is called
the {\it Hamming graph} $H(d,n)$. The Hamming graph is
regular. The number of equidistant vertices is denoted by
\begin{equation}\label{deg}
  \kappa_i:=|\{x'\in \mathcal{X}:\rho(x,x')=i\}|
  =\binom{d}{i}(n-1)^i, \quad i\in\{0,1,\ldots,d\}
\end{equation}
for any $x\in\mathcal{X}$, where $\kappa_1$ is the degree
of a vertex, namely, the number of adjacent vertices.
The Hamming graph is a fundamental object in algebraic
combinatorics, and group-theoretic approach is considered
to be effective \cite{BBIT21}. In this paper, we consider
Gaussian free fields on Hamming graphs. The Hamming
graph is quite different from the integer lattice
$\mathbb{Z}^d$ with edges between nearest-neighbour
vertices. For example, the diameter of $H(d,n)$ is $d$,
while that of a cube in $\mathbb{Z}^d$ with side length
$n$ is $(n-1)d$, despite having the same number of
vertices $n^d$.
The two graphs have $(n-1)dn^d/2$ and $(n-1)dn^{d-1}$
edges respectively.
We can thus say that the Hamming graph is more densely
connected than the integer lattice, while being
``smaller'' in terms of diameter.

Let $\partial\subset \mathcal{X}$ be a distinguished set
of vertices, called the {\it boundary} of the graph,
and set $\mathcal{Y}=\mathcal{X}\setminus\partial$.
In this paper, we discuss a class of discrete Gaussian
free fields on $H(d,n)$, denoted by
$(g_x)_{x\in\mathcal{X}}\in\mathbb{R}^{\mathcal X}$,
defined via the {\it Hamiltonian}
for a given configuration
$(\xi_x)_{x\in\mathcal{X}}\in\mathbb{R}^\mathcal{X}$:
\begin{equation}\label{ham}
  H_\mathcal{Y}\{(\xi_x)_{x\in\mathcal{X}}\}
  :=\frac{1}{4}\sum_{i=1}^dJ_i\sum_{(x,x')\in R_i}
  (\xi_x-\xi_{x'})^2
  +\frac{m^2}{2}\sum_{x\in \mathcal{Y}}\xi_x^2,
  \quad \mathcal{Y} \subset \mathcal{X},
\end{equation}
where $(\xi_x)_{x\in\mathcal{X}}$
represents a realisation of the Gaussian
free field $(g_x)_{x\in\mathcal{X}}$, referred to as
the {\it spin} configuration. The Hamiltonian tends to
be small if spins take similar values. Here, we assume that
the interaction between spins $g_x$ and $g_{x'}$
depends solely on the Hamming distance
$i=\rho(x,x')$ through weights $J_i\ge 0$, which satisfy
$\sum_{i=1}^d \kappa_i J_i\le 1$. The parameter $m\ge 0$
is called {\it mass}. For $\beta>0$, a centered Gaussian
free field ($\mathbf{E}[g_x]=0$, $x\in\mathcal{X}$) 
that satisfies the boundary condition $g_x=0$ for
$x\in\partial$ is given by the Gibbs measure:
\begin{equation}\label{gibbs}
  \mu_\mathcal{Y}(S)=\frac{1}{Z_\mathcal{Y}}\int_S
  e^{-\beta H_\mathcal{Y}\{(\xi_x)_{x\in\mathcal{X}}\}}
  \prod_{x\in\mathcal{Y}} d\xi_x
  \quad \text{for~each} \quad S\in\mathcal{B},
\end{equation}
where $\mathcal{B}$ is the Borel set over the hyperplane
\begin{equation}\label{plane}
  \mathcal{P}_\mathcal{Y}=\{(\xi_x)_{x\in\mathcal{X}}\subset
  \mathbb{R}^\mathcal{X}:\xi_x=0~\text{for~all}~x\in\partial\}.
\end{equation}
The normalisation constant $Z_\mathcal{Y}$ is called
the {\it partition function.} In studies of spin systems,
$\beta$ is called the {\it inverse temperature.}
The effect of this scaling is trivial, since \eqref{ham}
is a homogeneous polynomial; more precisely, under this
scaling, the variance of the Gaussian free field becomes
$1/\beta$ times
that of the unscaled field. Unless otherwise specified,
the interaction weights $J_i$, $i\in\{1,\ldots,d\}$ are
assumed not to depend on $m$ or $\beta$.

Research on discrete Gaussian free fields has a long
history; however, much of this research has focused on
their properties on the integer lattice $\mathbb{Z}^d$,
particularly in two dimensions ($d=2$). In particular,
there has been considerable research in the context of
Liouville quantum gravity and Schramm--Loewner
evolution (see \cite{BP25} for a book-length treatment),
as well as in the context of extreme values
\cite{Biskup20}. On the other hand, relatively little
research has focused on Gaussian free fields on graphs
other than $\mathbb{Z}^d$. Examples include Gaussian
free fields on compact manifolds \cite{CG20} and
the discrete torus \cite{Abacherli19}. To the best of
the author's knowledge, Gaussian free fields on
Hamming graphs have not been studied. The hypercube
$H(d,2)$ is a special case of a Hamming graph, and
spectral expansions of Gaussian free fields on
the hypercube in terms of independent standard
Gaussian random variables have been studied in
\cite{CG25,Griffiths25}.
In this paper, we discuss discrete Gaussian free
fields on $H(d,n)$ with homogeneous interactions
determined solely by the Hamming distance.

In Section~\ref{sect:rw}, we introduce random walks
for studying the class of Gaussian free fields.
First, we introduce a class of random walks on
Hamming graphs, where the transition probability is
determined solely by the Hamming distance between
vertices, and identify its Green function with
the covariance of the corresponding Gaussian free
fields. We then discuss the {\it massive} ($m>0$)
cases without boundary ($\partial=\emptyset$), and
the {\it massless} ($m=0$) cases with a boundary.
Here, the boundary refers to the outside of a set
of vertices equidistant from the origin. For
the former, we derive an expansion of the Green
function in terms of an orthonormal system and
provide an expression for the Gaussian free field
in terms of complex Gaussian
random variables. For the latter, we provide
an expression for the variance of the field at
the origin.

In Section~\ref{sect:part}, we discuss general
properties of the class of Gaussian free fields using
the random walk representation introduced in
Section~\ref{sect:rw}, In particular, we discuss
the massive cases without boundary. We provide
explicit expressions for the partition function
and demonstrate the convergence of the free energy,
defined in Section~\ref{sect:part}, to the limits for
large $d$ or $n$. We show that in the massive limit
($m\to\infty$), the scaled field
$(mg_x)_{x\in\mathcal{X}}$ converges weakly to i.i.d.\
Gaussian random variables. We then discuss
the covariance of Gaussian free fields. We show that
in the massless limit ($m\to 0$), Gaussian random
variables diverge, while the scaled field
$(mg_x)_{x\in\mathcal{X}}$ becomes perfectly correlated.
We also show that, in the limits of large $d$ or $n$,
the field $(g_x)_{x\in\mathcal{X}}$ converges weakly to
i.i.d.\ Gaussian random variables. We also show that,
for large $n$, the covariance decays exponentially with
the Hamming distance, irrespective of the interaction
weights.

In Section~\ref{sect:model}, we discuss three specific
models in order to investigate properties that cannot
be fully captured by general results of
Section~\ref{sect:part}. Section~\ref{sect:model:ind}
considers a model in which the interaction weights
between all pairs of vertices are equal. We provide
explicit results for the massive cases with and without
boundary, and for the massless case with a boundary.
In Section~\ref{sect:model:nn}, we consider a model
in which only adjacent vertices interact with equal
weights, i.e., the nearest-neighbour model. We show that,
in the massless case without boundary, the covariance
decays exponentially with the Hamming distance for
large $d$ or $n$. Section~\ref{sect:model:binom}
considers a model in which the interaction weights
follow a binomial probability mass function in
the distance with parameter
$\gamma\in(0,1-1/n]$, where $d(1-1/n)$ is the typical
distance of the Hamming graph.
We show that convergence of the variance
as $d\to\infty$ is fastest when the weights are
constant, i.e., when the model reduces to the first model.
Furthermore, we show that coupling $\gamma$ to $m$
yields a double-scaling limit in which the scaled field
$(mg_x)_{x\in\mathcal{X}}$ has correlations that
are convex in Hamming distance and is exactly
the reciprocal of the number of vertices at that
distance from a given vertex, attaining its minimum
strictly in the interior of $\{0,1,\ldots,d\}$. This
shows that the covariance structure of the Gaussian free
field completely encodes the geometry of the Hamming graph.

\section{The random walk representation}
\label{sect:rw}

In this section, we introduce random walks which are
useful to analyse Gaussian free fields on Hamming
graphs. See
\cite[Section 1.1]{BP25}, \cite[Sections 8.4,5]{FV17},
and \cite[Section 13.1]{Georgii88} for the random walk
representation. See \cite[Chapter 1,2]{BBIT21} for
details about the Hamming graph $H(d,n)$ and related
notions. 

Consider a discrete-time Markov chain $\{X_t:t\ge 0\}$
starting from $X_0=x\in\mathcal{X}$. We assume
a homogeneity of transition probability matrix $P$.

\begin{assumption}\label{assump}
  For $(x,x') \in \mathcal{X}\times\mathcal{X}$ and
  $(y,y') \in \mathcal{X}\times\mathcal{X}$,
  $\rho(x,x')=\rho(y,y') \Rightarrow 
  (P)_{x,x'}=(P)_{y,y'}$.
\end{assumption}

The wreath product of the symmetric groups of orders
$n$ and $d$, $\mathfrak{S}_n\wr\mathfrak{S}_d$, acts
on $\mathcal{X}=\{0,1,\ldots,n-1\}^d$ as
\begin{equation}\label{wreath}
  (\tau_1,\ldots,\tau_d;\sigma)(x_1,\ldots,x_d)
  =\big(\tau_1(x_{\sigma^{-1}(1)}),\ldots,\tau_d(x_{\sigma^{-1}(d)})\big)
\end{equation}
for each element
$g=(\tau_1,\ldots,\tau_d;\sigma)\in\mathfrak{S}_n\wr\mathfrak{S}_d$,
where $\tau_i\in\mathfrak{S}_n$ for each $i\in\{1,\ldots,d\}$
and $\sigma\in\mathfrak{S}_d$. Assumption~\ref{assump} is
the invariance of the transition probability under this
action, namely, $(P)_{g(x),g(x')}=(P)_{x,x'}$.

Such a transition probability matrix can be
represented as
\begin{equation}\label{trans}
  P=\sum_{i=0}^d\frac{w_i}{\kappa_i}A_i, \quad 
  w_i\ge 0, \quad \sum_{i=0}^d w_i=1,
\end{equation}
where $A_i$ is the {\it adjacency matrix} defined
by
\[
  (A_i)_{x,x'}=\left\{
    \begin{array}{ccc}
      1 & \text{if} & (x,x')\in R_i,\\ 
      0 & \text{if} & (x,x')\notin R_i,
    \end{array}\right.
\]
and $\kappa_i$ is defined in \eqref{deg}.
Since $P$ is symmetric, the uniform distribution is
a stationary distribution, and if $P$ is irreducible,
it is the unique stationary distribution.

Let $\hat{P}$ denote the restriction of $P$ to
$\mathcal{Y}\times\mathcal{Y}$, that is,
\[
    (\hat{P})_{x,x'}=\left\{
    \begin{array}{ll}     
       (P)_{x,x'} & \text{if} ~~ (x,x')
       \in\mathcal{Y}\times\mathcal{Y}\\
       0        & \text{otherwise},
    \end{array}
    \right.
\]
and $\tau_\partial$ denotes the first time that
the random walk $\{X_t:t\ge 0\}$ starting from
$X_0=x$ hits the boundary $\partial$, i.e.,
$\tau_\partial=\inf\{t\ge 0:X_t\in\partial\}$.
The Green function with killing is the expected
number of visits to $x'$ before hitting $\partial$
or being killed. Let the killing rate be $1-\alpha$,
$\alpha\in(0,1]$ and set 
$\tau_\Delta=\inf\{t\ge 1:X_t=\Delta\}$, where
$\Delta$ is the cemetery state. We have
\begin{equation}
  G_{\alpha,\partial}(x,x')
  =\mathbf{E}_x\Big(\sum_{t=0}^\infty
  1_{\{X_t=x', \tau_\partial>t, \tau_\Delta>t\}}\Big)
  =\sum_{t=0}^\infty \alpha^t(\hat{P}^t)_{x,x'}
  =(I-\alpha\hat{P})^{-1}_{x,x'} \label{green}
\end{equation}
for $(x,x')\in \mathcal{Y}\times\mathcal{Y}$, where
$I$ is the identity matrix. The last equality is
justified if the largest eigenvalue modulus of
$\alpha\hat{P}$ is strictly smaller than one.
Therefore, $\alpha\in(0,1)$ or $\partial\neq\emptyset$
is required. In fact, if $\alpha=1$ and 
$\partial=\emptyset$, we cannot define a Gaussian
free field (See Section~\ref{sect:part}).

When $\alpha=1$, i.e., there is no killing,
we must consider the case $\partial\neq\emptyset$.
While it is generally not easy to calculate
\eqref{green} exactly when the boundary is
not empty, it is relatively easy to calculate
for the case where $x'=x$. Set
$\tau_x=\inf\{t\ge 1:X_t=x\}$. Since the number
of returns to $x$ before hitting $\partial$
follows the geometric distribution of mean
$\mathbf{P}_x(\tau_\partial\le\tau_x)$, we have
\begin{equation}\label{green1}
  G_{1,\partial}(x,x)
  =\sum_{t=0}^\infty(\hat{P}^t)_{x,x}
  =\frac{1}{\mathbf{P}_x(\tau_\partial\le \tau_x)}.
\end{equation}

The following proposition is a version of
\cite[Proposition 1.8]{BP25} with the killing.
The proof is straightforward, but it is shown
for convenience.

\begin{proposition}
  Given $\mathcal{Y}\subset \mathcal{X}$ and
  boundary condition $g_x=0$, $x\in\partial$,
  let
  $(\xi_x)_{x\in\mathcal{X}}\in\mathbb{R}^\mathcal{X}$
  denote a realisation of the Gaussian
  free field $(g_x)_{x\in\mathcal{X}}$. Then,
  \begin{equation*}
    \frac{1}{4}\sum_{i=1}^d\frac{w_i}{\kappa_i}
    \sum_{(x,x')\in R_i}(\xi_x-\xi_{x'})^2
    +\frac{1-\alpha}{2\alpha}\sum_{x\in\mathcal{Y}}\xi_x^2
    =\frac{1}{2\alpha}\sum_{x,x'\in\mathcal{Y}}
    G_{\alpha,\partial}^{(-1)}(x,x')\xi_x\xi_{x'}.
  \end{equation*}
\end{proposition}

\begin{proof}
\begin{align*}  
  \sum_{x,x'\in\mathcal{Y}}
  &(1-\alpha\hat{P})_{x,x'}\xi_x\xi_{x'}
  =\sum_{x,x'\in\mathcal{X}}(1-\alpha P)_{x,x'}\xi_x\xi_{x'}\\
  =&-\frac{1}{2}\sum_{x,x'\in\mathcal{X}}
  (1-\alpha P)_{x,x'}(\xi_x-\xi_{x'})^2\\
  &+\frac{1}{2}\sum_{x,x'\in\mathcal{X}}
  (1-\alpha P)_{x,x'}\xi_x^2
  +\frac{1}{2}\sum_{x,x'\in\mathcal{X}}
  (1-\alpha P)_{x,x'}\xi_{x'}^2\\
  =&-\frac{1}{2}\sum_{x,x'\in\mathcal{X}}
  (1-\alpha P)(\xi_x-\xi_{x'})^2
  +\frac{\alpha}{2}\sum_{x,x'\in\mathcal{X}}
  (1-P)_{x,x'}\xi_x^2\\
  &+\frac{\alpha}{2}\sum_{x,x'\in\mathcal{X}}
  (1-P)_{x,x'}\xi_{x'}^2
  +\frac{1-\alpha}{2}\sum_{x,x'\in\mathcal{X}}
  (\delta_{x,x'}\xi_x^2+\delta_{x,x'}\xi_{x'}^2)\\
  =&\frac{\alpha}{2}\sum_{x,x'\in\mathcal{X}}
  (P)_{x,x'}(\xi_x-\xi_{x'})^2
  +(1-\alpha)\sum_{x\in\mathcal{Y}}\xi_x^2\\
  =&\frac{\alpha}{2}\sum_{i=1}^d
  \frac{w_i}{\kappa_i}\sum_{(x,x')\in R_i}
  (\xi_x-\xi_{x'})^2
  +(1-\alpha)\sum_{x\in\mathcal{Y}}\xi_x^2,
\end{align*}
where the second-to-last equality holds by
$\sum_{x'\in\mathcal{X}}(1-P)_{x,x'}=0$ for all
$x\in\mathcal{X}$ and $(P)_{x,x'}=(P)_{x',x}$.
Then, the assertion follows by \eqref{green}.
\end{proof}

%The first term of the left-hand side of \eqref{identity}
%is the Dirichlet energy.

The following proposition is a version of
\cite[Proposition 1.2]{BP25} with the killing.
A proof is similar to that of \cite{BP25}, so we
omit it.

\begin{proposition}\label{pro:nonneg}
  The function $G_{\alpha,\partial}$ is symmetric and
  non-negative definite. $G_{\alpha,\partial}(x,\cdot)$
  is discrete $\alpha$-harmonic in $\mathcal{Y}\setminus\{x\}$;
  more precisely $G_{\alpha,\partial}$ is the unique
  function such that
  $(I-\alpha\hat{P})G_{\alpha,\partial}(x,\cdot)=\delta_x(\cdot)$
  for all $x\in\mathcal{Y}$, and satisfies
  $G_{\alpha,\partial}(x,\cdot)=0$ on $\partial$ for all
  $x\in \mathcal{X}$.
\end{proposition}  

%The Hamiltonian is minimized by $\alpha$-harmonic
%functions. The Gaussian free field (GFF) can be viewed as
%a ``Gaussian perturbation of an $\alpha$-harmonic function.''
%More precisely, we have the following theorem.

%\begin{Theorem}[Markov property, cf. Theorem 1.10, BP2025]
%  Fix $U\Subset\mathcal{V}_{q,d}$. The Gaussian free
%  field $h_x$, $x\in\mathcal{V}_{q,d}$ can be decomposed
%  as
%  \[
%    h=g+\varphi,
%  \]
%  where $g$ is Gaussian free field on $U$ satisfying
%  $g_x=0$, $x\in\partial$, and $\varphi$ is $\alpha$-harmonic
%  in $U$. Moreover, $g_x$ and $\varphi$ are independent.
%\end{Theorem}

Now, identifying $J_i$, $i\in\{1,\ldots,d\}$ and
$m^2$ with $w_i/\kappa_i$, $i\in\{1,\ldots,d\}$ and
$(1-\alpha)/\alpha$ in \eqref{gibbs}, respectively,
we see that the Gibbs measure defines a centred Gaussian
free field $(g_x)_{x\in \mathcal{X}}$:
\[
  \mu_\mathcal{Y}(S)=\frac{1}{Z_\mathcal{Y}}
  \int_S\exp\Big\{-\frac{1}{2}
  \sum_{x,x'\in\mathcal{Y}}{\rm Cov}^{(-1)}(g_x,g_{x'})
  \xi_x\xi_{x'}\Big\}\prod_{x\in\mathcal{Y}} d\xi_x,
\]
where the covariance is identified as
${\rm Cov}(g_x,g_{x'})=\alpha G_{\alpha,\partial}(x,x')/\beta$.
By Proposition~\ref{pro:nonneg}, the covariance is
non-negative. We have
$w_0=1-\sum_{i=1}^dw_i=1-\sum_{i=1}^d\kappa_iJ_i$.
The identification allows us to investigate structure
of a Gaussian free field by analysing the associated
random walk. The Hamiltonian \eqref{ham} is expressed as
\begin{equation}\label{ham_rw}
  H_\mathcal{Y}((\xi_x)_{x\in\mathcal{X}})
  =\frac{1+m^2}{2}\sum_{x\in\mathcal{Y}}\xi_x^2
  -\frac{1}{2}\sum_{i=1}^d\frac{w_i}{\kappa_i}\sum_{(x,x')\in R_i}
  \xi_x\xi_{x'}
\end{equation}
on the hyperplane $\mathcal{P}_{\mathcal{Y}}$ given
in \eqref{plane}.

When a random field $(g_x)_{x\in \mathcal{X}}$ is
discrete $\alpha$-harmonic in $\mathcal{Y}$ and
satisfies a boundary condition $g_x=\varphi_x$,
$x\in\partial$ for a fixed $\varphi$, we say
$(g_x)_{x\in\mathcal{X}}$ solves
the $\alpha$-{\it Dirichlet problem} in $\mathcal{Y}$
with the $\alpha$-{\it discrete Laplacian}
$P-I/\alpha$. The following proposition is a version
of \cite[Theorem 8.26]{FV17} with the killing.

\begin{proposition}\label{prop:solve}
  The solution $(g_x)_{x\in\mathcal{X}}$ of
  the $\alpha$-Dirichlet problem with the boundary
  condition $g_x=\varphi_x$, $x\in\partial$, satisfies
  $\mathbf{E}[g_x]=\mathbf{E}_x[\varphi_{X_{\tau_\partial}}]$,
  $x\in \mathcal{X}$,
  where $\tau_{\partial}=\inf\{t\ge 0:X_t\in\partial\}$.
\end{proposition}

\begin{proof}
  If $x\in\partial$,
  $\mathbf{E}_x[\varphi_{X_{\tau_\partial}}]=\varphi_x$.
  Otherwise,
  \[
    \mathbf{E}_x[\varphi_{X_{\tau_\partial}}]=
    \sum_{y\in\mathcal{Y}}
    \mathbf{E}_x[\varphi_{X_{\tau_\partial}}|X_1=y]
    \mathbf{P}_x(X_1=y)
    =\sum_{y\in \mathcal{Y}}\alpha\hat{P}_{x,y}
    \mathbf{E}_y[\varphi_{X_{\tau_\partial}}],
  \]
  where the second equality holds by the strong Markov property.
  Since
  $\sum_{y\in\mathcal{Y}}(I-\alpha \hat{P})_{x,y}\mathbf{E}_y
  [\varphi_{X_{\tau_\partial}}]=0$, $x\in \mathcal{Y}$,
  $\mathbf{E}_x[\varphi_{X_{\tau_\partial}}]$
  is discrete $\alpha$-harmonic in $\mathcal{Y}$.
\end{proof}

Proposition~\ref{prop:solve} says that a centred Gaussian
free field defined by the Gibbs measure \eqref{gibbs}
solves the $\alpha$-Dirichlet problem with the boundary
condition $g_x=0$, $x\in\partial$.

In the remaining part of this section, we discuss
the following two settings:
\begin{itemize}
\item[(a)] $\alpha\in(0,1)$, corresponding to the massive
  case ($m>0$),
  with no boundary, $\partial=\emptyset$;
\item[(b)] $\alpha=1$, corresponding to massless case
  ($m=0$), with a boundary
    \begin{equation}\label{boundary}
      \partial_r
      =\{x\in \mathcal{X}:r+1\le \rho(x) \le d\},
      \quad r \in \{1,\ldots,d-1\},
    \end{equation}
\end{itemize}
where $\rho(x)=\rho(0,x)$. Without loss of generality,
we take the origin as the reference vertex. Then,
the set $\mathcal{Y}_r=\mathcal{X}\setminus\partial_r$
can be regarded as a ``ball'' of radius $r$ centred
at the origin.

\subsection{The massive case without boundary}

Hora \cite{Hor97} discussed the class of random walks
on $H(d,n)$ and gave a detailed treatment of the cut-off
phenomenon of a simple random walk ($w_i=\delta_{i,1}$).
The author and coauthor \cite{GM25a} discussed quantum
walks on $H(d,n)$, where the Fourier transform was
effectively used.

Let $P_t(i)$, ($t=0,1,\ldots$), denote the $t$-step
transition probability $(P^t)_{x,x'}$ for $(x,x')\in R_i$.
Hora \cite{Hor97} established a spectral representation
of the transition probability matrix
\begin{equation}\label{spec}
  P_t(j)=\frac{1}{n^d}\sum_{i=0}^d\lambda_i^t K_i(j),
  \quad j\in\{0,1,\ldots,d\}
\end{equation}
with eigenvalues
\begin{equation}\label{eigen} 
  \lambda_i=\sum_{j=0}^d \frac{w_j}{\kappa_j} K_j(i)
  \quad i\in\{0,1,\ldots,d\},
\end{equation}
where $\lambda_0=1$ and $-1/(n-1)\le \lambda_i\le 1$,
$i\in\{1,\ldots,d\}$. Here, 
\begin{equation}\label{Kpol}
  K_i(j)=\sum_{l=\max(0,i+j-d)}^{\min(i,j)}(-1)^l(n-1)^{i-l}
  \begin{pmatrix}j\\l\end{pmatrix}
  \begin{pmatrix}d-j\\i-l\end{pmatrix}
\end{equation}
are the Krawtchouk polynomials defined by the generating
function \cite[Section~2.10]{BBIT21}
\begin{equation}\label{gen}
  \sum_{l=0}^dK_l(j)s^l=(1+(n-1)s)^{d-j}(1-s)^j.
\end{equation}
In particular, $\kappa_i=K_i(0)$ and $K_0(j)=1$.
The lower bound for the eigenvalues follows from the
representation $\lambda_i=\mathbf{E}[\{-1/(n-1)\}^{L_i}]$,
$i\in\{1,\ldots,d\}$, where
\[
  \mathbf{P}(L_i=l)=\sum_{j=l}^{d+l-i}
  \frac{\binom{i}{l}\binom{d-i}{j-l}}{\binom{d}{j}}w_j,
  \quad l\in\{0,1,\ldots,i\}.
\]  
Here, $L_i$ follows a mixture of the hypergeometric
distributions.
Moreover, we have a useful bound
\begin{equation}\label{Kbound}
  \frac{|K_i(j)|}{\kappa_i}\le
  \sum_{l=\max(0,i+j-1)}^{\min(i,j)}(n-1)^{-l}
  \frac{\binom{j}{l}\binom{d-j}{i-l}}{\binom{d}{i}}
  =\mathbf{E}[(n-1)^{-{L_{i,j}}}]\le 1,
\end{equation}
where $L_{i,j}$ follows a hypergeometric distribution.
The orthogonality is with respect to
a binomial probability function:
\[
  \sum_{l=0}^dK_i(l)K_j(l)\binom{d}{l}
  \left(\frac{1}{n}\right)^{d-l}\left(1-\frac{1}{n}\right)^l
  =\delta_{i,j}\binom{d}{i}(n-1)^i.
\]
The spectral representation \eqref{spec} was obtained
by regarding the random walks as those on $\mathbb{Z}^d_n$
generated by circulant matrices \cite[Remark 27]{GM25b}.

We prepare a useful lemma that provides
the Fourier transform of a function of the Hamming
distance from the origin in terms of the Krawtchouk
polynomials. The assumption that $n$ is prime in
\cite[Lemma 4.6]{GM25a} is not needed here.

\begin{lemma}[\rm\cite{GM25a}]\label{lem:fourier}
  For any function $f:\{0,1,\ldots,d\}\to\mathbb{C}$,
  we have
  \begin{equation}\label{sum}
    \sum_{y\in \mathcal{X}}\zeta^{x\cdot y}f(\rho(y))
    =\sum_{i=0}^df(i)K_i(\rho(x)), \quad x\in \mathcal{X},
  \end{equation}
  where $x \cdot y=x_1 y_1+\cdots+x_d y_d$, and
  the $n$-th primitive root of unity is denoted by
  $\zeta\equiv e^{2\pi\sqrt{-1}/n}$.
\end{lemma}

If we choose $f(i)=w_i/\kappa_i$, the right-hand
side of \eqref{sum} is the same as that of
\eqref{eigen}, which means that
\[
  \lambda_{\rho(x)}=\mathbf{E}[\zeta^{x\cdot Y}],
  \quad \text{where} \quad
  \mathbf{P}(Y=y)=\frac{w_{\rho(y)}}{\kappa_{\rho(y)}},
  \quad y\in\mathcal{X}.
\]
%We note that $\zeta^{x\cdot y}$ is the character
%of the direct product of the cyclic groups
%$\mathbb{Z}_n^d$.

If $\alpha\in(0,1)$ ($m>0$) and $\partial=\emptyset$,
we have an expression of the Green function:
\begin{align}
  G_{\alpha,\emptyset}(x,x')
  &=\sum_{t=0}^\infty \alpha^t P_t(\rho(x,x'))
  =\frac{1}{n^d}\sum_{t=0}^\infty\sum_{j=0}^d
  (\alpha\lambda_j)^t K_j(\rho(x,x'))\nonumber\\
  &=\frac{1}{n^d}\sum_{j=0}^d\frac{K_j(\rho(x,x'))}
  {1-\alpha\lambda_j}. \label{green_h}
\end{align}

We introduce complex Gaussian random variables.
In the following, for $y\in\mathcal{X}$, we write
$n-y=(n-y_1,\ldots,n-y_d)$, where each component
is understood modulo $n$, so that is identified
with $0$.

\begin{definition}\label{def:hermite}
  Let $(Z_y)_{y\in\mathcal{X}}$ be a set of
  complex Gaussian random variables.
  For $y\in\mathcal{X}\setminus\{0\}$, the real and
  imaginary parts of $Z_y$ are independent centred
  Gaussian random variables with variance $1/2$,
  whereas $Z_0$ is real and follows the standard
  Gaussian distribution. In addition, they
  satisfy the reflection symmetry 
  ${\rm Re}Z_y={\rm Re}Z_{n-y}$ and skew-symmetry
  ${\rm Im}Z_y=-{\rm Im}Z_{n-y}$; equivalently,
  $Z_y=\overline{Z_{n-y}}$ for $y\in\mathcal{X}$.
  Apart from these symmetry relations,
  the corresponding real and imaginary parts are
  mutually independent. These random variables are
  referred to as the
  {\it Hermitian Gaussian random variables}
  on $\mathcal{X}$.
\end{definition}

\begin{proposition}\label{prop:hermite}
  The Hermitian Gaussian random variables $(Z_y)_{y\in\mathcal{X}}$
  on $\mathcal{X}$ satisfy
  $\mathbf{E}[Z_y\overline{Z_{y'}}]=\delta_{y,y'}$.
\end{proposition}

\begin{proof}
  We observe that
  $\mathbf{E}[Z_y\overline{Z_{y'}}]=
  \mathbf{E}[{\rm Re}Z_y{\rm Re}Z_{y'}+{\rm Im}Z_y{\rm Im}Z_{y'}]$.
  If $y=y'$, this is 1. Otherwise, if $y'\neq n-y$, then
  $\mathbf{E}[Z_y\overline{Z_{y'}}]=0$, since the centred
  random variables ${\rm Re}Z_y$, ${\rm Re}Z_{y'}$,
  ${\rm Im}Z_y$, and ${\rm Im}Z_{y'}$ are independent.
  If $y'= n-y$,
  \[
  \mathbf{E}[Z_y\overline{Z_{n-y}}]=\mathbf{E}[Z_y^2]=
  \mathbf{E}[({\rm Re}Z_y)^2-({\rm Im}Z_y)^2]+
  2\sqrt{-1}\mathbf{E}[({\rm Re}Z_y)({\rm Im}Z_y)]=0.
  \]  
\end{proof}

Now, we have the following representation theorem for
a Gaussian free field $(g_x)_{x\in\mathcal{X}}$.
This is an explicit version of \cite[Theorem~1.7]{BP25},
but involves Hermitian Gaussian random
variables on $\mathcal{X}$.

\begin{theorem}\label{thm:exp}
  Set $\alpha\in(0,1)$ $(m>0)$ and $\partial=\emptyset$.
  For $(x,x')\in \mathcal{X}\times\mathcal{X}$, we have
  the expansion of the Green function
  \begin{equation}\label{thm:exp1}
    \alpha G_{\alpha,\emptyset}(x,x')=
    \sum_{y\in\mathcal{X}}
    \frac{1}{1/\alpha-\lambda_{\rho(y)}}
    \frac{\zeta^{x\cdot y}}{\sqrt{n^d}}
    \frac{\zeta^{-x'\cdot y}}{\sqrt{n^d}}.
  \end{equation}
  Let $(Z_x)_{x\in \mathcal{X}}$ be
  Hermitian Gaussian random variables on $\mathcal{X}$.
  Then,
  \begin{equation}\label{thm:exp2}
    g_x=\frac{1}{\sqrt{\beta}}\sum_{y\in \mathcal{X}}
    \frac{Z_{y}}{\sqrt{1/\alpha-\lambda_{\rho(y)}}}
    \frac{\zeta^{x\cdot y}}{\sqrt{n^d}},
    \quad x\in\mathcal{X}
  \end{equation}
  is a real Gaussian free field given
  by the Gibbs measure \eqref{gibbs}. Here,
  $1/\alpha-\lambda_i$, $i\in\{0,1,\ldots,d\}$
  are the eigenvalues of $I/\alpha-P$
  $($the $\alpha$-discrete Laplacian with the sign
  reversed$)$ with an orthonormal basis of eigenfunctions
  $(\zeta^{x\cdot y}/\sqrt{n^d})_{y\in \mathcal{X}}$.
  The eigenvalues are degenerate; more precisely,
  for each $y\in\mathcal{X}$, the eigenfunction
  $\mathcal{X}\ni x\mapsto (\zeta^{x\cdot y}/\sqrt{n^d})\in\mathbb{C}$
  corresponds to the eigenvalue
  $1/\alpha-\lambda_{\rho(y)}$ whose multiplicity is
  $\kappa_{\rho(y)}$.
\end{theorem}

\begin{proof}
  Let $z_i\equiv(x_i-x'_i)\mod n$ for
  $i\in\{1,\ldots,d\}$. We observe
  $\rho(z)=\rho(x,x')$ and 
  $\zeta^{(x-x')\cdot y}=\zeta^{z\cdot y}$.
  By using the Fourier transform \eqref{sum},
  the expression of the Green function \eqref{green_h}
  can be recast to
  \[
    \alpha G_{\alpha,\emptyset}(x,x')
    =\frac{1}{n^d}\sum_{i=0}^d\frac{K_i(\rho(z))}
    {1/\alpha-\lambda_i}
    =\frac{1}{n^d}\sum_{y\in X}
    \frac{\zeta^{z\cdot y}}{1/\alpha-\lambda_{\rho(y)}}
    =\frac{1}{n^d}\sum_{y\in X}
    \frac{\zeta^{(x-x')\cdot y}}{1/\alpha-\lambda_{\rho(y)}}.
  \]
  For the second assertion, setting \eqref{thm:exp2},
  we confirm that $g_x$ is real:
  \[
    \overline{g_x}=\frac{1}{\sqrt{\beta}}\sum_{y\in\mathcal{X}}
    \frac{\overline{Z_y}}{\sqrt{1/\alpha-\lambda_{\rho(y)}}}
    \frac{\zeta^{-x\cdot y}}{\sqrt{n^d}}=
        \frac{1}{\sqrt{\beta}}\sum_{y\in\mathcal{X}}
    \frac{Z_{n-y}}{\sqrt{1/\alpha-\lambda_{\rho(n-y)}}}
    \frac{\zeta^{x\cdot (n-y)}}{\sqrt{n^d}}=g_x,
  \]
  where we used $\rho(y)=\partial(0,y)=\partial(n-y,0)=\partial(0,n-y)=\rho(n-y)$.
  Then, we have
  \begin{align*}
  {\rm Cov}(g_x,g_{x'})&=\mathbf{E}[g_x\overline{g_{x'}}]
  =\frac{1}{\beta n^d}\sum_{y\in\mathcal{X}}\sum_{y'\in\mathcal{X}}
  \frac{\mathbf{E}[Z_y\overline{Z_{y'}}]\zeta^{x\cdot y-x'\cdot y'}}
       {\sqrt{1/\alpha-\lambda_{\rho(y)}}
         \sqrt{1/\alpha-\lambda_{\rho(y')}}}\\
       &=
  \frac{1}{\beta n^d}\sum_{y\in\mathcal{X}}
  \frac{\zeta^{(x-x')\cdot y}}
       {1/\alpha-\lambda_{\rho(y)}}=\frac{\alpha}{\beta} G_{\alpha,\emptyset}(x,x'),      
  \end{align*}
  where we used Proposition~\ref{prop:hermite} in the
  third equality.
  The orthonormality of
  $(\zeta^{x\cdot y}/\sqrt{n^d})_{y\in \mathcal{X}}$
  is followed by
  \[
  \sum_{x\in \mathcal{X}} \zeta^{x\cdot (y-y')}
  =\sum_{x\in \mathcal{X}} \prod_{i=1}^d \zeta^{x_i (y_i-y'_i)}
  =\prod_{i=1}^d
  \Big(\sum_{x_i=0}^{n-1}\zeta^{x_i (y_i-y'_i)}\Big)
  =\prod_{i=1}^d (n\delta_{y_i,y'_i})=n^d\delta_{y,y'}.
  \]
  Using \eqref{sum} again, we expand the eigenvalue
  \eqref{eigen} as
  \begin{equation}\label{thm:exp3}
    \lambda_{\rho(x)}=\sum_{i=0}\frac{w_i}{\kappa_i}
    K_i(\rho(x))=  
    \sum_{y\in \mathcal{X}}\zeta^{x\cdot y}
    \frac{w_{\rho(y)}}{\kappa_{\rho(y)}}
  \end{equation}
  and the last assertion is followed by
  \[
    \sum_{x'\in \mathcal{X}}(P)_{x,x'}\zeta^{x'\cdot y}
    =\sum_{x'\in \mathcal{X}}\frac{w_{\rho(x,x')}}{\kappa_{\rho(x,x')}}
      \zeta^{x'\cdot y}
    =\sum_{z\in \mathcal{X}}\frac{w_{\rho(z)}}{\kappa_{\rho(z)}}\zeta^{(x-z)\cdot y}
    =\lambda_{\rho(y)}\zeta^{x \cdot y},
  \]
  where the complex conjugate of \eqref{thm:exp3} is
  used in the last equality. Since
  $\sum_{x'\in \mathcal{X}}(I/\alpha-P)_{x,x'}\zeta^{x'\cdot y}/\sqrt{n^d}=(1/\alpha-\lambda_{\rho(y)})\zeta^{x \cdot y}/\sqrt{n^d}$,
  we conclude that
  $(\zeta^{x\cdot y}/\sqrt{n^d})_{y\in\mathcal{X}}$,
  is an orthonormal basis of eigenfunctions. For each
  $y\in\mathcal{X}$, the eigenfunction
  $\mathcal{X}\ni x\mapsto (\zeta^{x\cdot y}/\sqrt{n^d})\in\mathbb{C}$
  corresponds to the eigenvalue
  $1/\alpha-\lambda_{\rho(y)}$ whose multiplicity is
  $\kappa_{\rho(y)}=|\{x\in\mathcal{X}:\rho(x)=\rho(y)\}|$. 
\end{proof}

\subsection{The massless case with a boundary}
\label{sect:rw:m0}

For the case of $m=0$ ($\alpha=1$), we consider
the boundary $\partial_r$, $r \in \{1,\ldots,d-1\}$
defined by \eqref{boundary}.
%In section~\ref{sect:model},
%the Green function $G_{1,\partial_r}$ will be used to
%discuss covariances of Gaussian random variables.

A gambler's ruin type argument gives the following lemma.

\begin{lemma}\label{lem:var}
  Set $\alpha=1$ $(m=0)$ and the boundary
  $\partial=\partial_r=\{x\in \mathcal{X}:r+1\le \rho(x) \le d\}$,
  $r \in \{1,\ldots,d-1\}$,
  and $h_i=\{x\in\mathcal{X}:\rho(x)=i\}$,
  $i\in\{0,1,\ldots,d\}$. The Green function is given as
  \begin{equation}\label{var_0}
    G_{1,\partial_r}(0,0)=\frac{1}{\sum_{i=1}^r w_ip_i+\sum_{i=r+1}^dw_i},
  \end{equation}
  where $(p)_{i\in\{1,\ldots,r\}}$ solves
  the system of linear equations
  \begin{equation}\label{lem:var:1}
    (I-\tilde{P})p=\tilde{p}_\partial, 
  \end{equation}
  where 
  $(\tilde{P})_{i,j\in\{1,\ldots,r\}}$ are the transition
  probabilities of the lumped chain of \eqref{trans}
  with the Hamming distance, whose state space is
  $h_i$, $i\in\{0,1,\ldots,d\}$, and 
  $(\tilde{p}_\partial)_i=1-\sum_{j=0}^r(\tilde{P})_{ij}$.
\end{lemma}

\begin{proof}
  The expression of the Green function \eqref{green1}
  gives
  $G_{1,\partial_r}(0,0)=1/\mathbf{P}_0(\tau_{\partial_r}\le\sigma_0)$, where $\tau_{\partial_r}=\inf\{t\ge 0:X_t=\partial_r\}$
  and $\sigma_0=\inf\{t\ge 1:X_t=0\}$.
  We can compute $\mathbf{P}_0(\tau_{\partial_r}\le\sigma_0)$
  as follows. The first step should be a jump to
  either of $x\in\mathcal{Y}\setminus\{0\}$ or
  hit $\partial_r$. By using the strong Markov property,
  we have
  \begin{align}
    \mathbf{P}_0(\tau_{\partial_r}\le \sigma_0)
    &=\sum_{x\in\mathcal{Y}\setminus\{0\}} (P)_{0,x}
    \mathbf{P}_x(\tau_{\partial_r}\le \sigma_0)
    +\sum_{x\in\partial_r} (P)_{0,x}\nonumber\\
    &=\sum_{x\in\mathcal{Y}\setminus\{0\}} (P)_{0,x}
    \mathbf{P}_x(\tau_{\partial_r}\le \sigma_0)+\sum_{i=r+1}^d w_i.
    \label{lem:var:2}
  \end{align}
  We denote the equivalence classes induced by
  the Hamming distance from the origin as
  $h_i=\{x\in\mathcal{X}:\rho(x)=i\}$,
  $i\in\{0,1,\ldots,d\}$.
  These sets are invariant under the action of
  the wreath product of symmetric groups
  $\mathfrak{S}_{n-1}\wr\mathfrak{S}_d$
  introduced by \eqref{wreath}, but here
  $\mathfrak{S}_{n-1}$ is the symmetric group of
  $\{1,2,\ldots,n-1\}$.
  The lumping of the transition probability matrix
  \eqref{trans} exists if 
  \begin{equation*}%\label{lem:var:3}
    \sum_{x'\in h_j} (P)_{x,x'}=\sum_{x'\in h_j} (P)_{y,x'}
  \end{equation*}
  is satisfied for each $x,y\in h_i$. This equation
  holds, because there exists
  $g\in\mathfrak{S}_{n-1}\wr\mathfrak{S}_d$
  such that $y=g(x)$, and 
  \[
    \sum_{x'\in h_j} (P)_{x,x'}
    =\sum_{x'\in h_j} (P)_{g(x),g(x')}
    =\sum_{x'\in h_j} (P)_{y,g(x')}
    =\sum_{x'\in h_j} (P)_{y,x'},
  \]
  where the first equality follows by \eqref{assump}.
  Therefore, the transition probabilities of the lumped
  chain are given by
  $(\tilde{P})_{\rho(x)j}=\sum_{x'\in h_j} (P)_{x,x'}$,
  $x\in\mathcal{X}$. Substituting
  $p_i=\mathbf{P}_x(\tau_{\partial_r}\le\sigma_0)$,
  $x\in h_i$, $i\in\{1,\ldots,r\}$ into
  \eqref{lem:var:2} gives
  $\mathbf{P}_0(\tau_{\partial_r}\le\sigma_0)
    =\sum_{i=1}^r w_i p_i+\sum_{i=r+1}^d w_i$.
  and \eqref{var_0} follows. By the strong Markov property,
  we have
  \begin{equation*}
    \mathbf{P}_0(\tau_{\partial_r}>\sigma_0)=1-p_i=
    \sum_{j=1}^r (\tilde{P})_{ij}(1-p_j)+(\tilde{P})_{i0},\quad x\in h_i,
  \end{equation*}
  and this is the system \eqref{lem:var:1}.
\end{proof}

We have a simple proposition.

\begin{proposition}\label{prop:ind}
  If the transition probabilities of a random walk
  satisfying Assumption~\ref{assump} are independent of
  the current state, namely, the condition
  \begin{equation}\label{cor:ind:cond}
    (P)_{x,x'}=(P)_{y,x'} \quad \text{for~any}
    \quad (x,y) \in \mathcal{X}\times\mathcal{X},
  \end{equation}
  is satisfied for each $x'\in\mathcal{X}$,
  the transition probabilities are
  $(P)_{x,x'}=1/n_d$ for all $x,x'\in\mathcal{X}$.
  Moreover, we have
  \begin{equation}\label{cor:ind:G}
    G_{1,\partial_r}(0,0)
    =\frac{1+|\partial_r|}{|\partial_r|}.
  \end{equation}
\end{proposition}

\begin{proof}
  Since the transition probability matrix is symmetric,
  we have
  $(P)_{x',x}=(P)_{x'y}$ for any
  $(x,y) \in \mathcal{X}\times\mathcal{X}$,
  which means that all entries of the transition matrix 
  are $1/n_d$. Then, $(\tilde{P})_{ij}=w_j=\kappa_j/n^d$ and
  the system \eqref{lem:var:1} reduces to the equations:
  $p_i=\sum_{j=1}^rw_jp_j+\sum_{j=r+1}^dw_j$,
  $i \in\{1,\ldots,r\}$.
  This is solved immediately:
  \[
    p_i=\frac{1-\sum_{j=0}^r w_j}{1-\sum_{j=1}^r w_j}
    =\frac{n^d-\sum_{i=0}^r \kappa_i}
    {n^d-\sum_{i=1}^r \kappa_i}=\frac{|\partial_r|}{1+|\partial_r|},
    \quad i\in\{1,\ldots,r\},
  \]    
  and we have the assertion.
\end{proof}

\section{The partition function and covariance}
\label{sect:part}

In this section, we present some general properties
of the class of Gaussian free fields on Hamming graphs
given by the Gibbs measure \eqref{gibbs}. We discuss
the massive case without boundary using the random walk
representation introduced in Section~\ref{sect:rw}.
The massive cases with a boundary are discussed in
Section~\ref{sect:model}, as their properties are
model-specific.

The following proposition shows that the Gaussian free
fields on Hamming graphs given by the Gibbs measure
\eqref{gibbs} do not exhibit spontaneous magnetisation,
i.e., a non-zero expectation of the spatial average.
This is expected since we consider centred fields.

\begin{proposition}\label{pro:magnet}
  Set $m>0$ and $\partial=\emptyset$. The spatial average
  of a Gaussian free field given by the Gibbs measure
  \eqref{gibbs} is a centred Gaussian random variable with
  variance $1/(n^d\beta m^2)$.
\end{proposition}

\begin{proof}
  By using the expansion \eqref{thm:exp1}, the mean of
  the Gaussian free field becomes
  \begin{align*}
    \frac{1}{n^d}\sum_{x\in\mathcal{X}}g_x
    &=\frac{1}{n^d\sqrt{\beta}}\sum_{x\in\mathcal{X}}\sum_{y\in\mathcal{X}}
    \frac{Z_y}{\sqrt{1/\alpha-\lambda_{\rho(y)}}}
    \frac{\zeta^{x\cdot y}}{\sqrt{n^d}}\\
    &=\frac{1}{n^d\sqrt{\beta}}\sum_{y\in\mathcal{X}}
    \frac{Z_y}{\sqrt{1/\alpha-\lambda_{\rho(y)}}}
    \frac{n^d\delta_{\rho(y),0}}{\sqrt{n^d}}
    =\frac{1}{\sqrt{n^d}\sqrt{\beta}}\frac{Z_0}{\sqrt{1/\alpha-1}},
  \end{align*}             
  where the third equality follows by \eqref{sum} with
  $f=1$ and \eqref{gen}. Since $1/\alpha=1+m^2$ and $Z_0$ is
  the standard Gaussian random variable, the assertion
  follows.
\end{proof}

\subsection{The partition function}

%A group-theoretic argument gives
We have a simple expression for the partition function.
If $m=0$ and $\partial=\emptyset$, the partition
function $Z_\mathcal{X}$ does not exist, so we
cannot define a Gaussian free field. We have
already seen in Section~\ref{sect:rw} that
we cannot define the Green function
$G_{1,\emptyset}(x,x')$.

\begin{proposition}\label{prop:part}
  Set $m>0$ and $\partial=\emptyset$. The partition
  function of a Gaussian free field given by
  the Gibbs measure \eqref{gibbs} is represented as
  \begin{equation}\label{prop:part:z}
    Z_\mathcal{X}=
    \left(\frac{2\pi}{\beta}\right)^{n^d/2}
    \frac{1}{m}\prod_{i=1}^d
    (1+m^2-\lambda_i)^{-\kappa_i/2},
  \end{equation}
  where $\lambda_i$ are the eigenvalues of
  the random walk \eqref{eigen} and $\kappa_i$
  is given by \eqref{deg}.
\end{proposition}

\begin{proof}%[Proof of Theorem \ref{thm:part}]
  The partition function is given by the Gaussian integral
  \begin{align*}
    Z_\mathcal{X}
    &=\int_{\mathbb{R}^\mathcal{X}}
    e^{-\beta H_\mathcal{X}\{(\xi_x)_{x\in\mathcal{X}}\}}\mathrm{d}\xi
    =\int_{\mathbb{R}^\mathcal{X}}
    e^{-\frac{\beta}{2\alpha} \sum_{x,x'\in\mathcal{X}}
      G^{(-1)}_{\alpha,\emptyset}(x,x')\xi_x\xi_{x'}}\mathrm{d}\xi\\
    &=\left(\frac{2\pi\alpha}{\beta}\right)^{n^d/2}
    |G_{\alpha,\emptyset}|^{1/2},
  \end{align*}
  where $|G_{\alpha,\emptyset}|$ is the determinant of
  the matrix whose $(x,x')$-entry is $G_{\alpha,\emptyset}(x,x')$.
  Since $\partial=\emptyset$, using the expression \eqref{green}
  with $\hat{P}=P$, we obtain
  \begin{equation}\label{prop:part:det}
  |G_{\alpha,\emptyset}|=|I-\alpha P|^{-1}
  =\prod_{y\in\mathcal{X}}\frac{1}{1-\alpha\lambda_{\rho(y)}}.
  \end{equation}
  Therefore, 
  \begin{equation*}
    Z_{\mathcal{X}}
    =\left(\frac{2\pi}{\beta}\right)^{n^d/2}\prod_{y\in\mathcal{X}}
    \frac{1}{\sqrt{1/\alpha-\lambda_{\rho(y)}}}
    =\left(\frac{2\pi}{\beta}\right)^{n^d/2}\frac{1}{m}
    \prod_{i=1}^d(1+m^2-\lambda_i)^{-\kappa_i/2}
  \end{equation*}
  where we used $\lambda_0=1$ and the multiplicities of
  the eigenvalues given in Theorem~\ref{thm:exp}.
\end{proof}

\begin{remark}\label{rem:thermo}
  The factor $(2\pi/\beta)^{n^d/2}$ is the partition
  function of the {\it Maxwell--Boltzmann distribution}
  describing particle speeds in idealised gases,
  where the particles move freely without
  interacting with one another. The distribution
  is that of independent Gaussian random variables
  of variance $1/\beta$. The {\it internal energy}
  in statistical mechanics is obtained as
  $-\partial\log {Z}_\mathcal{X}/{\partial\beta}={n^d}/(2\beta)$,
  which coincides with that of the Maxwell--Boltzmann
  distribution, where $n^d$ corresponds to the number
  of particles. In this sense, the determinant
  \eqref{prop:part:det} encodes the interaction
  among particles. The (Helmholtz) free energy is
  defined as $(-\log Z_{\mathcal{X}})/\beta$.
\end{remark}

\begin{corollary}\label{cor:part}
  Set $m>0$ and $\partial=\emptyset$.
  A Gaussian free field given by the Gibbs measure
  \eqref{gibbs} satisfies the following.
  \begin{itemize}
  \item[$(i)$] If $m\to\infty$, the scaled field
    $(\sqrt{\beta}mg_x)_{x\in\mathcal{X}}$ converges
    weakly to a collection of independent standard
    Gaussian random
    variables.
  \item[$(ii)$] Assume the stochastic matrix
    \eqref{trans} is irreducible. As $m\to 0$,
    \[
      mZ_\mathcal{X}\to \left(\frac{2\pi}{\beta}\right)^{n^d/2}
      \prod_{i=1}^d(1-\lambda_i)^{-\kappa_i/2}.
    \]
  \end{itemize}
\end{corollary}

\begin{proof}
  (i) By the expression \eqref{prop:part:z}, we
  observe $Z_\mathcal{X}\sim\{2\pi/(\beta m^2)\}^{n^d/2}$
  as $m\to\infty$, because the eigenvalues
    $\lambda_i$, $i\in\{1,\ldots,d\}$ are fixed.
  The asymptotic form is independent Gaussian random
  variables of variance $1/(\beta m^2)$. Since
  the partition function characterizes the distribution,
  the assertion holds. (ii) It is obvious from
  \eqref{prop:part:z}.
\end{proof}

The limit distribution as $m\to\infty$ is known
and is consistent with Proposition~\ref{pro:magnet},
whereas the limit as $m\to 0$ is non-trivial.
We are therefore interested in the limit as $m\to 0$.
In the next section, we also discuss the case $m=0$.
Since we cannot define a massless Gaussian free field
without boundary, we incorporate the boundary $\partial_r$
as in Section~\ref{sect:rw:m0}.

Next, we consider other limits, where the number
of vertices $n^d$ approaches infinity. Note that in
both limits the number of vertices grows to infinity,
but not only does it grow at a different rate,
but also the diameter of the Hamming graph is given
by $d$, so in the limit as $n\to\infty$ the diameter
does not change.

\begin{remark}\label{rem:thermo_lim}
  Both the limits $d\to\infty$ and $n\to\infty$
  can be called {\it thermodynamic limits,}
  because the number of particles grows to
  infinity (see Remark~\ref{rem:thermo}).
  For a Gaussian free field on an infinite set such
  as $\mathbb{Z}^d$, we take the limit of extending
  a finite subset with the outside of the finite
  subset as the boundary. Since the Hamming graph is
  a finite set, in the next section we discuss
  the limits as $d\to\infty$ or $n\to\infty$ without
  boundary, as well as the joint limit as $r,d\to\infty$
  with the boundary $\partial_r$.
  %Note that we cannot take $r\to\infty$ while fixing $d$.  
\end{remark}

\begin{corollary}\label{cor:lim}
  Set $m>0$ and $\partial=\emptyset$. The free energy of
  a Gaussian free field given by the Gibbs measure
  \eqref{gibbs} per particle satisfies the following.
  \begin{itemize}
  \item[$(i)$]     
  Assume there exists $\lambda_*=\lim_{d\to\infty}\lambda_{\lfloor(1-1/n)d\rfloor}$, and that, for every $\epsilon>0$, there exists
  $\delta>0$ such that
  $\max_{i\in A_d(\epsilon)}|\lambda_i-\lambda_*|\to 0$ as
  $d\to\infty$, where
  $A_d(\epsilon)=\{i\in\{0,1,\ldots,d\}:|i-(1-1/n)d|<\epsilon d^{1/2+\delta}\}$.
  Then,
  \begin{equation*}
    -\frac{1}{\beta}
    \lim_{d\to\infty}\frac{\log Z_{\mathcal{X}}}{n^d}
    =\frac{1}{2\beta}\left\{\log(1+m^2-\lambda_{*})
    -\log\left(\frac{2\pi}{\beta}\right)\right\};
  \end{equation*}
  \item[$(ii)$]
  Assume there exists
  $\lambda_{\#}=\lim_{n\to\infty}\lambda_d$. Then,
  \begin{equation*}
    -\frac{1}{\beta}
    \lim_{n\to\infty}\frac{\log Z_{\mathcal{X}}}{n^d}
    =\frac{1}{2\beta}\left\{\log(1+m^2-\lambda_\#)-
    \log\left(\frac{2\pi}{\beta}\right)\right\}.
  \end{equation*}
  \end{itemize}
\end{corollary}  

\begin{proof}
  The expression \eqref{prop:part:z} leads to
  \begin{equation}\label{col:lim1}
    -\frac{\log Z_\mathcal{X}}{n^d/2}=
    \sum_{i=0}^d{d\choose i}
    \Big(1-\frac{1}{n}\Big)^i
    \Big(\frac{1}{n}\Big)^{d-i}
    \log(1+m^2-\lambda_i)-\log\left(\frac{2\pi}{\beta}\right)
  \end{equation}
  and the limit of the summation is our concern.\\
  (i) The first term of the right-hand side of
  \eqref{col:lim1} is
  $\mathbf{E}[\log(1+m^2-\lambda_{I_d})]$, where
  $I_d$ follows a binomial distribution. Since
  $-1/(n-1)\le\lambda_i\le 1$, we have a bound
  $\log m^2 \le \log(1+m^2-\lambda_i) \le \log\{1+m^2+1/(n-1)\}$
  for all $i\in\{0,1,\ldots,d\}$. Then,
  \begin{align*}
  &\mathbf{E}[|\log(1+m^2-\lambda_{I_d})-\log(1+m^2-\lambda_*)|]\\
  &\le
  \max_{i\in A_d(\epsilon)}|\log(1+m^2-\lambda_i)-\log(1+m^2-\lambda_*)|+
  C\mathbf{P}(I_d\in A^c_d(\epsilon)),
  \end{align*}
  where $C=\log[1+\{1+1/(n-1)\}/m^2]$. By using 
  $\max_{i\in A_d(\epsilon)}|\log(1+m^2-\lambda_i)-\log(1+m^2-\lambda_*)|\to 0$
  as $d\to\infty$ and
  $\mathbf{P}(I_d\in A^c_d(\epsilon))\le e^{-2\epsilon^2d^{2\delta}}$
  by Hoeffding's inequality, we have
  $\lim_{d\to\infty}\mathbf{E}[\log(1+m^2-\lambda_{I_d})]=\log(1+m^2-\lambda_*)$.
  (ii) It is obvious from the expression.
\end{proof}
  
\subsection{The covariance}

We are interested in the decay of covariances of
Gaussian free fields with the Hamming distance.
If $\partial=\emptyset$, we saw that the limit
distribution of a scaled field $(mg_x)_{x\in\mathcal{X}}$
as $m\to\infty$ is a collection of i.i.d.\ 
Gaussian random variables. But the limiting distribution
as $m\to 0$ seems non-trivial
(see Corollary~\ref{cor:part}, ii). 
The following corollary gives more detailed information
in this respect.

\begin{corollary}\label{cor:cov}
  Set $m>0$ and $\partial=\emptyset$ and assume
  the stochastic matrix \eqref{trans} is
  irreducible. If $m\to 0$, a Gaussian free field
  $(g_x)_{x\in\mathcal{X}}$ given by the Gibbs
  measure \eqref{gibbs} satisfies
  $n^dm^2 {\rm Cov}(g_x,g_{x'})\to 1/\beta$.
  In this limit, the covariances of the Gaussian
  random variables diverge, while the scaled field
  $(mg_x)_{x\in\mathcal{X}}$ is perfectly correlated.
\end{corollary}

\begin{proof}
  By the Perron--Frobenius theorem, $\lambda_0=1$
  and $\lambda_i<1$, $i\in\{1,\ldots,d\}$. From
  the expansion \eqref{thm:exp1}, we have
  \[
  m^2n^d{\rm Cov}(g_x,g_{x'})=
  m^2n^d\frac{\alpha}{\beta}G_{\alpha,\emptyset}(x,x')=
  \frac{1}{\beta}\left\{1
    +\sum_{y\in \mathcal{X}\setminus\{0\}}
    \frac{m^2\zeta^{(x-x')\cdot y}}
    {1-\lambda_{\rho(y)}+m^2}\right\}.
  \]
  The second term satisfies
  \[
    \left|\sum_{y\in\mathcal{X}\setminus\{0\}}
    \frac{m^2\zeta^{(x-x')\cdot y}}{1-\lambda_{\rho(y)}+m^2}\right|
    =\sum_{y\in\mathcal{X}\setminus\{0\}}
    \frac{m^2}{1-\lambda_{\rho(y)}+m^2}
    \le \frac{m^2(n^d-1)}{\epsilon+m^2}\to 0
  \]
  as $m^2n^d\to 0$, where $\epsilon>0$ is one minus
  the second largest eigenvalue, where we note that
  we assumed that $n^d\ge 4$. Therefore,
  $m^2n^d{\rm Cov}(g_x,g_{x'})\to 1/\beta$ and
  ${\rm Corr}(mg_x,mg_{x'})\to 1$.
\end{proof}

\begin{remark}\label{rem:deloc}
  Corollary~\ref{cor:cov} suggests that for any
  interaction satisfying Assumption~\ref{assump},
  the covariance diverges as $m\to 0$.
  This contrasts sharply with Gaussian free fields
  on the integer lattice $\mathbb{Z}^d$ with
  nearest-neighbour interactions.
  In the thermodynamic limit of such models
  (see Remark~\ref{rem:thermo_lim}),
  the covariance decays exponentially with
  the Euclidean distance for any $d\ge 1$
  \cite[Section 8.5]{FV17}. In this respect,
  further analysis is provided in
  Section~\ref{sect:model:nn:m}.
  Note also that we assume here that the interaction
  weights are independent of $m$. In fact, as we shall
  see in Section~\ref{sect:model:binom}, if
  the interaction weights are scaled with $m$,
  the claims of Corollary~\ref{cor:cov} do not hold, and
  a non-trivial correlation structure emerges in
  the limit.
\end{remark}

The next proposition concerns the asymptotic behaviour
as $n\to\infty$. It states that, for any interaction
satisfying Assumption~\ref{assump}, the covariance
decays exponentially with the Hamming distance,
provided that $C_{m,d}(\rho(x,x'))>0$.

\begin{proposition}\label{pro:limcov}
  Set $m>0$ and $\partial=\emptyset$. For
  $(x,x')\in\mathcal{X}\times\mathcal{X}$,
  the covariance of a Gaussian free field
  $(g_x)_{x\in\mathcal{X}}$ given by the Gibbs measure
  \eqref{gibbs} satisfies
  \begin{equation}\label{pro:limcov:1}
    {\rm Cov}(g_x,g_{x'})
    =\frac{C_{m,d}(\rho(x,x'))}{\beta}
     e^{-{\rho(x,x')\log n}}
        +O(n^{-\rho-1}), \quad n\to\infty,
  \end{equation}
  where
  \[
    C_{m,d}(\rho)=\lim_{n\to\infty}
    \sum_{i=0}^\rho
    \binom{\rho}{i}\frac{(-1)^i}{1+m^2-\lambda_{i+d-\rho}}.
  \]
\end{proposition}

\begin{proof}
  The expression of the Krawtchouk polynomials
  \eqref{Kpol} gives
  \[
  K_i(j)=\left\{
  \begin{array}{ll}
      (n-1)^i\binom{d-j}{i}\{1+O(n^{-1})\} & j\le d-i,\\
      (-1)^{i+j-d}(n-1)^{d-j}\binom{j}{i+j-d}
      \{1+O(n^{-1})\} & j>d-i. 
   \end{array} \right.  
  \]
  Substituting this into \eqref{green_h}, we have
  \begin{align*}
    \alpha G_{\alpha,\emptyset}(x,x')&=
    \frac{(n-1)^{d-\rho}}{n^d}\left\{
    \frac{1}{1/\alpha-\lambda_{d-\rho}}
    +\sum_{i=d-\rho+1}^d
    \binom{\rho}{i+\rho-d}
    \frac{(-1)^{i+\rho-d}}
         {1/\alpha-\lambda_i}\right\}\\
    &     +O(n^{-\rho-1})\\
    &=C_{m,d}(\rho)n^{-\rho}+O(n^{-\rho-1}),
  \end{align*}
  where $\rho=\rho(x,x')$, and the assertion
  holds.
\end{proof}

The asymptotic behaviour as $d\to\infty$ is more subtle
because it requires a uniform evaluation of the summand
in the Green function \eqref{green_h}. To discuss this,
we need to specify the interaction weights, as we will
see in Section~\ref{sect:model}. However, the following
explicit result can be obtained for the limit.
  
\begin{proposition}\label{pro:limcovd}
  Let $\lambda_*$ be the constant defined in Corollary~\ref{cor:lim}.
  Set $m>0$ and $\partial=\emptyset$.
  For $(x,x')\in\mathcal{X}\times\mathcal{X}$, the covariance
  of a Gaussian free field $(g_x)_{x\in\mathcal{X}}$ given
  by the Gibbs measure \eqref{gibbs} satisfies
  \begin{equation*}
    \lim_{d\to\infty}{\rm Cov}(g_x,g_{x'})=\frac{1}{\beta}
    \frac{\delta_{x,x'}}{1+m^2-\lambda_*}.
  \end{equation*}
\end{proposition}

\begin{proof}
  If $x=x'$, the Green function \eqref{green_h} gives
  \begin{align*}
    {\rm Var}(g_x)
    &=\frac{\alpha}{\beta}G_{\alpha,\emptyset}(x,x)
    =\frac{1}{\beta n^d}\sum_{i=0}^d\frac{K_i(0)}
    {1/\alpha-\lambda_i}
    =\frac{1}{\beta n^d}\sum_{i=0}^d\frac{\kappa_i}
    {1/\alpha-\lambda_i}\\
    &=\frac{1}{\beta}\sum_{i=0}^d
    {d \choose i}\left(1-\frac{1}{n}\right)^i
    \left(\frac{1}{n}\right)^{d-i}
    \frac{1}{1/\alpha-\lambda_i}.
  \end{align*}
  The same argument used in the proof of
  Corollary~\ref{cor:lim} gives
  \[
  \lim_{d\to\infty}{\rm Var}(g_x)
  =\frac{1}{\beta}\frac{1}{1+m^2-\lambda_*}.
  \]
  If $x\neq x'$, we have
  \begin{equation*}
    \alpha G_{\alpha,\emptyset}(x,x')=\frac{1}{n^d}
    \sum_{i=0}^d \frac{K_i(\rho)}{1/\alpha-\lambda_i}
    =\sum_{i=0}^d
    \binom{d}{i}\left(1-\frac{1}{n}\right)^i
    \left(\frac{1}{n}\right)^{d-i}
    \frac{K_i(\rho)/\kappa_i}{1+m^2-\lambda_i},
  \end{equation*}
  where $\rho=\rho(x,x')$.
  Substituting $s=1$ into \eqref{gen}, we have $\sum_{i=0}^d K_i(\rho)=0$
  for $\rho>0$. Therefore, we may write
  \begin{align*}
    \alpha G_{\alpha,\emptyset}(x,x')
    &=\sum_{i=0}^d
    \binom{d}{i}\left(1-\frac{1}{n}\right)^i
    \left(\frac{1}{n}\right)^{d-i}
    \frac{K_i(\rho)}{\kappa_i}
    \left(\frac{1}{1+m^2-\lambda_i}-\frac{1}{1+m^2-\lambda_*}\right)\\
    &=\mathbf{E}\left[
    \frac{K_{I_d}(\rho)}{\kappa_{I_d}}
    \left(\frac{1}{1+m^2-\lambda_{I_d}}-\frac{1}{1+m^2-\lambda_*}
    \right)\right],
  \end{align*}
  where $I_d$ follows a binomial distribution. Then,
  \begin{align*}
    |\alpha G_{\alpha,\emptyset}(x,x')|&\le
    \mathbf{E}\left[\left|
    \frac{1}{1+m^2-\lambda_{I_d}}-\frac{1}{1+m^2-\lambda_*}
    \right|\right]\\
    &\le
    \max_{i\in A_d(\epsilon)}
    \left|\frac{1}{1+m^2-\lambda_{I_d}}-\frac{1}{1+m^2-\lambda_*}\right|
    +C'\mathbf{P}(I_d\in A_d^c(\epsilon))
  \end{align*}
  for a positive constant $C'$ independent of $d$, where we used
  \eqref{Kbound} in the first inequality, and the bound
  $-1/(n-1)\le \lambda_i \le 1$, $i\in\{0,1,\ldots,d\}$ in
  the second inequality. By using
  $\max_{i\in A_d(\epsilon)}|1/(1+m^2-\lambda_{I_d})-1/(1+m^2-\lambda_*)|\to 0$
  as $d\to\infty$ and
  $\mathbf{P}(I_d\in A_d^c(\epsilon))\le e^{-2\epsilon^2 d^{2\delta}}$
  by Hoeffding's inequality, we have
  $\lim_{d\to\infty}|\alpha G_{\alpha,\emptyset}(x,x')|=0$,
  where $A_d(\epsilon)$ is defined in Corollary~\ref{cor:lim}.
  Therefore, we conclude that
  $\lim_{d\to\infty}{\rm Cov}(g_x,g_{x'})=\lim_{d\to\infty}\alpha G_{\alpha,\emptyset}(x,x')/\beta=0$ for $x\neq x'$.
\end{proof}

This result is in accordance with Corollary~\ref{cor:lim}
(i) for the limit of the free energy, because it
coincides with that of the independent Gaussian random
variables of variance $1/\{\beta(1+m^2-\lambda_*)\}$.

As we saw in Corollary~\ref{cor:cov}, the covariance
of the Gaussian random variables diverges as $m\to 0$.
However, after taking the limit $d\to\infty$ or
$n\to\infty$, the variance remains finite, even as $m\to 0$.
Moreover, Propositions~\ref{pro:limcov} and
\ref{pro:limcovd} can be summarised in the following theorem.

\begin{theorem}\label{thm:limcov}
  Let $\lambda_*$ be the constant defined in
  Corollary~\ref{cor:lim} and assume there exists
  $\lambda_\flat=\lim_{n\to\infty}\lambda_{d-\rho(x,x')}$.
  Set $m>0$ and $\partial=\emptyset$.
  For $(x,x')\in\mathcal{X}\times\mathcal{X}$,
  the covariance of a Gaussian free field given by
  the Gibbs measure \eqref{gibbs} satisfies 
  \begin{equation*}
    \lim_{d\to\infty} {\rm Cov}(g_x,g_{x'})=
    \frac{1}{\beta}\frac{\delta_{x,x'}}{1+m^2-\lambda_*},
    \quad
    \lim_{n\to\infty} {\rm Cov}(g_x,g_{x'})=
    \frac{1}{\beta}\frac{\delta_{x,x'}}{1+m^2-\lambda_\flat}.
  \end{equation*}  
  For both limits, the Gaussian free field converges weakly
  to a collection of i.i.d.\ Gaussian random variables.
\end{theorem}

In summary, when $m>0$ and $\partial=\emptyset$,
we have the following. As $m\to 0$, the covariance of
a Gaussian free field diverges, and the scaled field
$(mg_x)_{x\in\mathcal{X}}$ becomes perfectly correlated
(Corollary~\ref{cor:cov}), whereas the scaled field
converges weakly to i.i.d.\ Gaussian random variables
as $m\to\infty$ (Corollary~\ref{cor:part}, i).
However, after taking the limit $d\to\infty$ or
$n\to\infty$, the Gaussian free field converges weakly
to a collection of i.i.d.\ Gaussian random variables,
even as $m\to 0$. Moreover, the covariance decays
exponentially with the Hamming distance as $n\to\infty$
(Proposition~\ref{pro:limcov}). These properties hold
regardless of the interaction weights, provided that
they satisfy Assumption~\ref{assump}.

\section{Specific models}
\label{sect:model}

In this section, we discuss three specific models of
Gaussian free fields on Hamming graphs. We explicitly
present the results obtained in Sections~\ref{sect:rw}
and \ref{sect:part} and investigate further properties
of these models.

\subsection{The distance-independent interaction}
\label{sect:model:ind}

We consider the Hamiltonian
\begin{equation}\label{ham:mean}
  H_\mathcal{Y}\{(\xi_x)_{x\in\mathcal{X}}\}
  =\frac{1}{4n^d}\sum_{x,x'\in\mathcal{X}}
  (\xi_x-\xi_{x'})^2
  +\frac{m^2}{2}\sum_{x\in\mathcal{Y}} \xi_x^2,
  \quad \mathcal{Y}\subset\mathcal{X},
\end{equation}
in which every pair of vertices interacts with equal
weight $J_i=1/n^d$, $i\in\{1,\ldots,d\}$. This is
a mean-field model analogous to the Curie--Weiss model
in discrete spin systems (see \cite[Chapter 2]{FV17}).

According to the recipe obtained in Section~\ref{sect:rw},
we consider the random walk determined by
$w_i/\kappa_i=1/n^d$, $i\in\{0,1,\ldots,d\}$,
and the eigenvalues \eqref{eigen} are
$\lambda_i=\delta_{i,0}$, $i\in\{0,1,\ldots,d\}$.
This random walk chooses the next vertex uniformly
from all vertices. Thus, it is equivalent to the simple
random walk on the complete graph with self-loops.

The following results are independent of the Hamming
distance between the vertices. Although some expressions
involve the boundary $\partial_r$, where $r$ is
the Hamming distance from the origin, these expressions
remain valid for any choice of boundary.

\subsubsection{The massive case without boundary}

If $m>0$ and $\partial=\emptyset$, the partition
function \eqref{prop:part:z} becomes
\[
  Z_\mathcal{X}=
  \left(\frac{2\pi}{\beta}\right)^{n^d/2}
  \frac{1}{m}(1+m^2)^{-(n^d-1)/2}
  =\left\{\frac{2\pi}{\beta(1+m^2)}\right\}^{n^d/2}
  \sqrt{\frac{1+m^2}{m^2}}.
\]
The free energy converges to the same limit as $d\to\infty$
or $n\to\infty$ (see Corollary~\ref{cor:lim}):
\[
  -\frac{1}{\beta}\frac{\log Z_{\mathcal{X}}}{n^d}\to
  \frac{1}{2\beta}\left\{\log(1+m^2)-
  \log\left(\frac{2\pi}{\beta}\right)\right\}.
\]
The covariance is given explicitly by the Green
function \eqref{green_h} as 
\begin{align*}
  {\rm Cov}(g_x,g_{x'})
  &=\frac{\alpha}{\beta}G_{\alpha,\emptyset}(x,x')
  =\frac{\alpha}{\beta n^d}\Big\{
  \frac{1}{1-\alpha}+\sum_{i=1}^d K_i
  (\rho(x,x'))\Big\}\\
  &=\frac{\alpha}{\beta n^d}\Big\{\frac{1}{1-\alpha}+(n^d\delta_{x,x'}-1)
  \Big\}\\
  &=\frac{1}{\beta}\left\{\frac{\delta_{x,x'}}{1+m^2}+\frac{1}{n^d}
  \frac{1}{(1+m^2)m^2}\right\}, \quad (x,x')\in\mathcal{X}\times\mathcal{X},
\end{align*}
where the second equality follows by \eqref{green_h},
and the third equality follows by \eqref{gen}.
This expression is consistent with
Theorem~\ref{thm:limcov}, but we note that
$C_{m,d}(\rho(x,x'))=0$ in Proposition~\ref{pro:limcov}
if $0<\rho(x,x')<d$.

\subsubsection{The massive case with a boundary}

Since the transition probabilities of the random walk
are independent of the current state, as
\eqref{cor:ind:cond}, we can obtain the Green function
even for the case with a boundary.

Let $\partial=\partial_r$, $r\in\{1,\ldots,d-1\}$.
For each step, there are the following four
possible events for the random walk at
$x\in\mathcal{Y}=\mathcal{X}\setminus\partial_r$.
We set $h_i=\{x\in\mathcal{X}:\rho(x)=i\}$.
\begin{itemize}
\item[(i)] visits $x'\in\mathcal{Y}$ with probability $\alpha/n^d$;
\item[(ii)] hits $\partial_r$ with probability
  $\alpha\sum_{i=r+1}^d w_i$;
\item[(iii)] visits
  $\mathcal{Y}\setminus\{x'\}=\mathcal{X}\setminus(\partial_r\cup\{x'\})$
  with probability
  $\alpha(1-1/n^d-\sum_{i=r+1}^dw_i)$;
\item[(iv)] is killed with probability $1-\alpha$.
\end{itemize}
If (ii) or (iv) occurs, the random walk stops.
Therefore, the expected number of visits $x'$
before stopping can be computed as
\begin{equation*}
  \delta_{x,x'}+
  \sum_{t=1}^\infty\sum_{k=0}^{t-1}
  k{t-1\choose k}\left(\frac{\alpha}{n_d}\right)^k
  \left\{\alpha
  \left(1-\frac{1}{n^d}-\sum_{i=r+1}^dw_i\right)\right\}^{t-k-1}
  \left\{1-\alpha+\alpha\sum_{i=r+1}^dw_i\right\}.
\end{equation*}
The summation is straightforward, and we obtain
\[
  G_{\alpha,\partial_r}(x,x')=\delta_{x,x'}
  +\frac{1}{n^d}\frac{1}{1/\alpha-1-|\partial_r|/n^d},
  \quad \alpha\in(0,1)
\]
and the covariance is
\begin{equation*}
  {\rm Cov}(g_x,g_{x'})
  =\frac{1}{\beta}\left\{\frac{\delta_{x,x'}}{1+m^2}
  +\frac{1}{n^d}\frac{1}{(1+m^2)(m^2+|\partial_r|/n^d)}
  \right\}, \quad (x,x')\in(\mathcal{X}\setminus\partial_r)
  \times(\mathcal{X}\setminus\partial_r).
\end{equation*}
%The assertions of Theorem~\ref{thm:limcov} still hold
%even when $\partial\neq\emptyset$. However,
In contrast to Corollary~\ref{cor:cov}, the limit as
$m\to 0$ remains finite. The following result is
an immediate consequence, although its behaviour is
counter-intuitive: the contribution to the covariance
from the mode that would be the zero mode if
$\partial=\emptyset$, whose eigenvalue is instead $m^2$
in the expansion \eqref{thm:exp2}, remains finite and
decreases as the boundary grows.

\begin{proposition}
  Set $m>0$ and $\partial=\partial_r$. For
  $(x,x')\in(\mathcal{X}\setminus\partial_r)\times(\mathcal{X}\setminus\partial_r)$, the Gaussian free
  field $(g_x)_{x\in\mathcal{X}}$ defined via the
  Hamiltonian \eqref{ham:mean} satisfies
  \begin{equation}\label{cov_indb}
  \lim_{m\to 0}{\rm Cov}(g_x,g_{x'})
  =\frac{1}{\beta}\left(\delta_{x,x'}+\frac{1}{|\partial_r|}\right).
  \end{equation}
  Consequently, for $x\neq x'$,
  $\lim_{m\to 0}{\rm Corr}(g_x,g_{x'})=1/(1+|\partial_r|)$,
  which is strictly decreasing in the boundary size $|\partial_r|$.
\end{proposition}

We can compute the partition function by direct
integration:
\begin{align}
  Z_\mathcal{Y}&=\int_{\mathcal{P}^\mathcal{Y}}
    e^{-\beta H_\mathcal{Y}((\xi_x)_{x\in\mathcal{X}})}
    \prod_{x\in\mathcal{Y}}\mathrm{d}\xi_x\nonumber\\
  &=\int_{\mathcal{P}^\mathcal{Y}}
    \exp\Big\{{-\frac{1+m^2}{2}\beta\sum_{x,x'\in\mathcal{Y}}
    (I-aJ)_{x,x'}\xi_x\xi_{x'}}\Big\}
    \prod_{x\in\mathcal{Y}}\mathrm{d}\xi_x\nonumber\\
  &=\int_{\mathcal{P}_\mathcal{Y}}\exp\Big\{-\frac{1+m^2}{2}\beta
    \big(\eta_1^2+\cdots+\eta_{|\mathcal{Y}|-1}^2
    +(1-|\mathcal{Y}|a)\eta_{|\mathcal{Y}|}^2\big)
    \Big\}\prod_{x\in\mathcal{Y}}\mathrm{d}\eta_x\nonumber\\
    &=\left(\frac{2\pi}{\beta(1+m^2)}\right)^{|\mathcal{Y}|/2}
    \sqrt{\frac{1+m^2}
      {|\partial_r|/n^d+m^2}},
    \quad \frac{1}{a}=(1+m^2)n^d, \label{part_indb}
\end{align}
where we used the expression \eqref{ham_rw}, and $J$
is the matrix of all entries one. The third equality
follows by the fact that a matrix $(\delta_{i,j}-aJ)$
is diagonalised by an orthogonal matrix as
$U {\rm diag}(1,1,\ldots,1,1-|\mathcal{Y}|a) U^{-1}$.

\subsubsection{The massless case with a boundary}

If $m=0$ and $\partial=\partial_r$, the Green
function $G_{1,\partial_r}(0,0)$ is already
given in \eqref{cor:ind:G}, but a more general
result can be obtained immediately from
\eqref{green}
\[
  G_{1,\partial_r}(x,x')=\delta_{x,x'}+
  \sum_{t=1}^\infty
  \left(1-\frac{|\partial_r|}{n^d}\right)^{t-1}
  \frac{1}{n^d}=
  \delta_{x,x'}+\frac{1}{|\partial_r|},
\]
and we have
${\rm Cov}(g_x,g_x')=G_{1,\partial_r}(x,x')/\beta
  =\left(\delta_{x,x'}+1/|\partial_r|\right)/\beta$.
Since this expression coincides with \eqref{cov_indb},
the covariance is continuous at $m=0$.
The partition function is given by \eqref{part_indb},
and we have
$Z_{\mathcal{Y}}
  =(2\pi/\beta)^{|\mathcal{Y}|/2}
  \sqrt{n^d/|\partial_r|}$.
    
\subsection{The nearest-neighbour interaction}
\label{sect:model:nn}

We consider the Hamiltonian
\begin{equation}\label{ham:simple}
  H_\mathcal{Y}\{(\xi_x)_{x\in\mathcal{X}}\}=\frac{1}{4}
  \sum_{(x,x')\in R_1}\frac{1}{(n-1)d}(\xi_x-\xi_{x'})^2
  +\frac{m^2}{2}\sum_{x\in\mathcal{Y}}\xi_x^2,
  \quad \mathcal{Y}\subset\mathcal{X},
\end{equation}
in which nearest-neighbour pairs of vertices interact with
equal weight $1/\{(n-1)d\}$. Since models on the integer
lattice with nearest-neighbour interactions are prevalent
in the literature, it is natural to consider this model.

We consider the random walk determined by
$w_i/\kappa_i=\delta_{i,1}/\{(n-1)d\}$,
$i\in\{0,1,\ldots,d\}$,
and the eigenvalues \eqref{eigen} are
$\lambda_i=1-ni/\{(n-1)d\}$, $i\in\{0,1,\ldots,d\}$.
This is the simple random walk in each step moves
to a nearest-neighbour vertex. 

\subsubsection{The massive case without a boundary}
\label{sect:model:nn:m}

If $m>0$ and $\partial=\emptyset$, the partition
function \eqref{prop:part:z} becomes
\[
  Z_\mathcal{X}=\left(\frac{2\pi}{\beta}\right)^{n^d/2}\frac{1}{m}
  \prod_{i=1}^d
  \left(m^2+\frac{ni}{(n-1)d}\right)^{-\kappa_i/2}
\]
and the free energy converges to the same limit as
$d\to\infty$ or $n\to\infty$:
\[
  -\frac{1}{\beta}\frac{\log Z_{\mathcal{X}}}{n^d}\to
  \frac{1}{2\beta}\left\{\log(1+m^2)-
  \log\left(\frac{2\pi}{\beta}\right)\right\},
\]
where the assumption of Corollary~\ref{cor:lim} is
satisfied, because
$\max_{i\in A_d(\epsilon)}|\lambda_i-\lambda_*|\le n\epsilon d^{\delta-1/2}/(n-1)\to 0$
as $d\to\infty$ for $\delta\in (0,1/2)$, where
$\lambda_*=\lim_{d\to\infty}\lambda_{\lfloor (1-1/n)d\rfloor}=0$.

The covariance is given by the Green function
\eqref{green_h} as
\begin{equation}\label{cov:simp}
  {\rm Cov}(g_x,g_{x'})=\frac{\alpha}{\beta}G_{\alpha,\emptyset}(x,x')=
  \frac{1}{\beta n^d}\sum_{i=0}^d
  \frac{K_i(\rho(x,x'))}{m^2+ni/\{(n-1)d\}}, \quad
  (x,x')\in\mathcal{X}\times\mathcal{X}.
\end{equation}
The covariance decays exponentially with the Hamming distance
as $n\to\infty$ (Proposition~\ref{pro:limcov}), and we
have \eqref{pro:limcov:1} with
\begin{align*}
  C_{m,d}(\rho)
  &=\sum_{i=0}^\rho\binom{\rho}{i}\frac{(-1)^i}{1+m^2+(i-\rho)/d}
  =\sum_{i=0}^\rho\binom{\rho}{i}(-1)^i
  \int_0^1v^{m^2+(i-\rho)/d}\mathrm{d}v\\
  &=\int_0^1 v^{m^2-\rho/d}(1-v^{1/d})^\rho \mathrm{d}v
  =d\int_0^1 z^{d(1+m^2)-\rho-1}(1-z)^\rho \mathrm{d}z\\
  &=d~{\rm Beta}(d(1+m^2)-\rho,\rho+1).
\end{align*}
The following result describes the limit as $d\to\infty$.

\begin{proposition}\label{prop:simp:lim}
  Set $m>0$ and $\partial=\emptyset$.
  For $(x,x')\in\mathcal{X}\times\mathcal{X}$
  the Gaussian free field $(g_x)_{x\in\mathcal{X}}$
  defined via the Hamiltonian \eqref{ham:simple}
  satisfies the following.
  \begin{itemize}
  \item[(i)]
  For a fixed Hamming distance $\rho=\rho(x,x')$, we have
  \[
    {\rm Cov}(g_x,g_{x'})
    =\frac{\rho!}{\beta\{1+m^2\}^{\rho+1}}
    e^{-\rho\log \{(n-1)d\}}\{1+O(d^{-1})\},
    \quad d\to\infty.
  \]
  \item[(ii)]
    For $\rho(x,x')=\sigma d$, $\sigma\in(0,1]$, the covariance
    ${\rm Cov}(g_x,g_{x'})$ decays exponentially in $\sigma$
    as $d\to\infty$. Moreover, the exponential decay rate is
    strictly increasing in $\sigma$.
  \end{itemize}
\end{proposition}

\begin{proof}
  (i) The expression \eqref{cov:simp} leads to
  \begin{align*}
    G_{\alpha,\emptyset}&(x,x')=\frac{1}{n^d}\sum_{i=0}^d
    \int_0^1 v^{-\alpha\big\{1-\frac{ni}{(n-1)d}\big\}}K_i(\rho)
    \mathrm{d}v
    =\frac{1}{n^d}\int_0^1
    v^{-\alpha}\sum_{i=0}^d
    \{v^{\frac{n\alpha}{(n-1)d}}\}^iK_i(\rho)\mathrm{d}v\\
    &=\frac{1}{n^d}
    \int_0^1v^{-\alpha}\big\{1+(n-1)v^{\frac{n\alpha}{(n-1)d}}\big\}^{d-\rho}
    \big\{1-v^{\frac{n\alpha}{(n-1)d}}\big\}^{\rho}\mathrm{d}v.
    %&=\frac{1}{n^d}\frac{(n-1)d}{n\alpha}
    %\int_0^1s^{\frac{(n-1)d}{n\alpha}(1-\alpha)-1}
    %\{1+(n-1)s\}^{d-\rho}(1-s)^{\rho}\mathrm{d}s,
  \end{align*}
  %We changed the integration variable
  %from $v$ to $s=v^{n\alpha/\{(n-1)d\}}$ and used \eqref{gen}
  %in the third equality. Since the dominant contribution to
  %the integral comes from around $s=1$,
  Changing the integration variable from $v$ to
  $u=d(1-v^{n\alpha/\{(n-1)d\}})$ yields
  \begin{equation}\label{prop:simp:lim1}
    \alpha G_{\alpha,\emptyset}(x,x')
    =\frac{1-1/n}{(nd)^{\rho}}
    \int_0^d \left(1-\frac{u}{d}\right)^{m^2d(1-1/n)-1}
    \left\{1-\Big(1-\frac{1}{n}\Big)\frac{u}{d}\right\}^{d-\rho}
    u^{\rho}
    \mathrm{d}u.
  \end{equation}
  For large $d$, the integral becomes
  \begin{align*}
    &\int_0^d \left(1-\frac{u}{d}\right)^{m^2d(1-1/n)-1}
    \left\{1-\Big(1-\frac{1}{n}\Big)\frac{u}{d}\right\}^{d-\rho}
    u^{\rho}\mathrm{d}u\\
    &=\int_0^\infty e^{-(1+m^2)(1-1/n)u}u^\rho
    \{1+O(d^{-1})\}\mathrm{d}u\\
    &=\rho!\left\{\frac{1}{(1+m^2)(1-1/n)}\right\}^{\rho+1}
    \{1+O(d^{-1})\}.
  \end{align*}
  Since
  ${\rm Cov}(g_x,g_{x'})=\alpha G_{\alpha,\emptyset}(x,x')/\beta$,
  the assertion holds.\\
  (ii) The integral in \eqref{prop:simp:lim1} is recast into
  \begin{align*}
    &d^{\rho(x,x')+1}\int_0^1\frac{e^{-df(s,\sigma)}}{1-s}\mathrm{d}s,
    \\
    &f(s,\sigma)=-m^2(1-1/n)\log(1-s)-(1-\sigma)\log\{1-(1-1/n)s\}
    -\sigma\log s,
  \end{align*}
  where $s=u/d$ and $\sigma=\rho(x,x')/d$. We have
  \[
  f'(s,\sigma)
  \equiv\frac{\partial f}{\partial s}
  =\frac{m^2(1-1/n)}{1-s}+\frac{(1-\sigma)(1-1/n)}{1-(1-1/n)s}
  -\frac{\sigma}{s}
  \]
  and
  \begin{align*}
    &f(s,\sigma)>0,\\
    &f''(s,\sigma)\equiv\frac{\partial^2 f}{\partial s^2}
  =\frac{m^2(1-1/n)}{(1-s)^2}+\frac{(1-\sigma)(1-1/n)^2}{1-(1-1/n)s}
  +\frac{\sigma}{s^2}>0, \quad (s,\sigma)\in(0,1)\times(0,1].
  \end{align*}
  Since $\lim_{s\to 0}f'(s,\sigma)=-\infty$ and
  $\lim_{s\to 1}f'(s,\sigma)=+\infty$, $f(s,\sigma)$ attains
  a unique minimum at, say, $\hat{s}_{\sigma}\in(0,1)$, for
  each $\sigma\in(0,1]$. The saddle-point approximation yields
  \begin{align*}
    &\alpha G_{\alpha,\emptyset}(x,x')=
    \frac{1-1/n}{1-\hat{s}_\sigma}
    \sqrt{\frac{2\pi d}{f''(\hat{s}_\sigma,\sigma)}}
    e^{-dR(\sigma)}\{1+O(d^{-1})\},\\
    &R(\sigma)=f(\hat{s}_\sigma,\sigma)+\sigma\log n>\sigma\log n.
  \end{align*}
  Since $\lim_{\sigma\to 0}\hat{s}_\sigma=0$, we obtain
  $\lim_{\sigma\to 0} R(\sigma)=0$. Moreover, by the envelope
  theorem, we have
  \begin{align*}
    \frac{\mathrm{d}R(\sigma)}{\mathrm{d}\sigma}
    &=\frac{\partial f}{\partial\sigma}+\log n
    =\log\{1-(1-1/n)\hat{s}_\sigma\}-\log \hat{s}_\sigma+\log n\\
    &=\log\left\{1+n\left(\frac{1}{\hat{s}_\sigma}-1\right)\right\}
    >0.
  \end{align*}
  Therefore,
  ${\rm Cov}(g_x,g_{x'})=\alpha G_{\alpha,\emptyset}(x,x')/\beta$
  decays exponentially in $\sigma$, and the exponential decay
  rate $R(\sigma)$ is strictly increasing in $\sigma$.
\end{proof}  

As we saw in Remark~\ref{rem:deloc}, the behaviour
of the Gaussian free field with nearest-neighbour
interaction on the Hamming graph $H(d,n)$ is quite
different from that on the integer lattice $\mathbb{Z}^d$. 
However, for large $d$ or $n$, the behavior of the Gaussian
free field on $H(d,n)$ becomes similar to that on
$\mathbb{Z}^d$, in the sense that its covariance decays
exponentially with the Hamming distance.

\subsubsection{The massive case with a boundary}

Here, we extend the definition \eqref{boundary} of
$\partial_r$ by setting $\partial_r=\emptyset$ for
$r\ge d$, so that $\tau_{\partial_r}\to\infty$ as
$r\to\infty$ for each fixed $d$. Since the simple
random walk satisfies
$\{\tau_{\partial_r}>t\}\subset\{\tau_{\partial_{r+1}}>t\}$,
the monotone convergence theorem and \eqref{green} yield,
for each fixed $d$,
\begin{align*}
  \lim_{r\to\infty}
  G_{\alpha,\partial_r}(x,x')&=
  \mathbf{E}_x\Big(\lim_{r\to\infty}
  \sum_{t=0}^\infty1_{\{X_t=x',\tau_{\partial_r}>t,\tau_{\Delta}>t\}}\Big)
  =\mathbf{E}_x\Big(
  \sum_{t=0}^\infty1_{\{X_t=x',\tau_{\Delta}>t\}}\Big)\\
  &=G_{\alpha,\emptyset}(x,x').
\end{align*}
Here, $\tau_{\partial_r}=\inf\{t\ge 0:X_t\in\partial_r\}$.
Since this holds for each $d$, we may then let $d\to\infty$;
Proposition~\ref{prop:simp:lim} gives
\begin{align}
  \lim_{d\to\infty}\lim_{r\to\infty}{\rm Cov}(g_x,g_{x'})
  &=\frac{\alpha}{\beta}\lim_{d\to\infty}
  \lim_{r\to\infty}G_{\alpha,\partial_r}(x,x')
  =\frac{\alpha}{\beta}\lim_{d\to\infty}
  G_{\alpha,\emptyset}(x,x')\nonumber\\
  &=\frac{\delta_{x,x'}}{\beta(1+m^2)}.\label{cov:simp_m}
\end{align}

\subsubsection{The massless case with a boundary}

If $m=0$ and $\partial=\partial_r$, $r\in\{1,\ldots,d-1\}$,
we use Lemma~\ref{lem:var}. The system of linear equations
\eqref{lem:var:1} is equivalent to the three-term recurrence
relation
\[
  p_{x+1}-p_x=\frac{x}{d-x}\frac{1}{n-1}(p_x-p_{x-1}),
  \quad \text{for} \quad x=1,2,\ldots,r
\]
with the boundary condition $p_0=0$, $p_{r+1}=1$.
It is straightforward to see
\[
  1=p_{r+1}=\sum_{i=0}^r {d-1 \choose i}^{-1}(n-1)^{-i}p_1,
\]
and we have
\[
  p_1^{-1}=\sum_{i=0}^r{d-1 \choose i}^{-1}(n-1)^{-i}.
\]
The expression \eqref{var_0} gives
\begin{equation}\label{var:simp_0}
  {\rm Var}(g_0)
  =\frac{1}{\beta}G_{1,\partial_r}(0,0)
  =\frac{1}{\beta}\frac{1}{p_1}
  =\frac{1}{\beta}\sum_{i=0}^r{d-1\choose i}^{-1}(n-1)^{-i},
\end{equation}
which strictly decreases with the boundary size, because
the boundary size $|\partial_r|$ strictly decreases with $r$.
We have the following limit.

%Evaluating the summation in \eqref{var:simp_0} seems not
%straightforward, but

\begin{proposition}\label{prop:cov:simp_0}
  Set $m=0$ and
  $\partial=\partial_r$, $r\in\{1,\ldots,d-1\}$.
  The Gaussian free field $(g_x)_{x\in\mathcal{X}}$
  defined via the Hamiltonian \eqref{ham:simple}
  satisfies
  \begin{equation}\label{prop:cov:simp_0:1}
    \lim_{d\to\infty}{\rm Var}(g_0)=\frac{1}{\beta}.
  \end{equation}
\end{proposition}

\begin{proof}
  Fix $x=r/d\in(0,1)$ and suppose $x<1-1/n$. Since
  \[
    {{d-1}\choose{i}} (n-1)^i \ge
    \left(\frac{(d-i)(n-1)}{i}\right)^i,
    \quad i\in\{0,\ldots,d-1\},
  \]
  we have
  \begin{align*}
    \frac{1}{\beta} &\le {\rm Var}(g_0)
    \le \frac{1}{\beta}\sum_{i=0}^r
    \left(\frac{i}{(d-i)(n-1)}\right)^i
    \le  
    \frac{1}{\beta}\sum_{i=0}^r
    \left(\frac{x}{(1-x)(n-1)}\right)^i\\
    &=\frac{1}{\beta}
    \frac{1-[x/\{(1-x)(n-1)\}]^{r+1}}
         {1-x/\{(1-x)(n-1)\}}
    \to\frac{1}{\beta}
  \end{align*}
  as $r\to\infty$, where the first inequality
  holds by the term $i=0$ of the summation
  in \eqref{var:simp_0}. Since \eqref{var:simp_0}
  is increasing with $r$, the limit holds for
  any $r<(1-1/n)d$. If $x\ge 1-1/n$, we divide
  the summation in \eqref{var:simp_0} into the two parts:
  \[
    S_1=\sum_{i=0}^{\lfloor (1-1/n)d\rfloor}
      {d-1\choose i}^{-1}(n-1)^{-i}, \quad
    S_2=\sum_{i=\lceil (1-1/n)d\rceil}^r
    {d-1\choose i}^{-1}(n-1)^{-i}.
  \]
  We have already seen that $\lim_{d\to\infty} S_1=1$.
  For the second part, noting that the summand increases,
  because
  \[
  n^{1-d}{d-1\choose i}(n-1)^i=
  \binom{d-1}{i}\left(1-\frac{1}{n}\right)^i
  \left(\frac{1}{n}\right)^{d-1-i}
  \]
  is the probability mass function of a binomial
  distribution, we have
  \[
    |S_2|<(r+1-\lceil(1-1/n)d\rceil)
    {d-1 \choose r}^{-1}(n-1)^{-r}.
  \]
  By using the Stirling formula, we have
  \[
    {d-1 \choose r}^{-1}(n-1)^{-r}\sim
    \sqrt{2\pi d}
    \exp\left[d\{x\log x+(1-x)\log(1-x)-x\log(n-1)\}\right]
  \]
  as $r,d\to\infty$. Since the exponent is
  always negative, $\lim_{d\to\infty}|S_2|=0$
  and the assertion holds.
\end{proof}

Since \eqref{cov:simp_m} with $x=x'$ converges to
\eqref{prop:cov:simp_0:1} as $m\to 0$, the variance
is continuous at $m=0$.

\subsection{The binomially-weighted interaction}
\label{sect:model:binom}

We consider the Hamiltonian
\begin{equation}\label{ham_bin}
  H_\mathcal{Y}\{(\xi_x)_{x\in\mathcal{X}}\}
  =\frac{1}{4}\sum_{i=1}^d{d\choose i}
  \frac{\gamma^i(1-\gamma)^{d-i}}{\kappa_i}
  \sum_{(x,x')\in R_i}(\xi_x-\xi_{x'})^2
  +\frac{m^2}{2}\sum_{x\in\mathcal{Y}}\xi_x^2,
  ~ \mathcal{Y}\subset\mathcal{X},
\end{equation}
in which interactions exist between every pair of
vertices, as in the distance-independent model
discussed in Section~\ref{sect:model:ind}, but with
distance-dependent weights. The weights are determined
by the binomial probability mass function with
parameter $\gamma\in{(0,1-1/n]}$, with
the weights of long-range interactions increasing
relative to those of short-range interactions as $\gamma$
increases.
Note that $d(1-1/n)$, the expectation of the binomial
distribution, is the typical distance in the Hamming
graph: the number of vertices at a given distance from
a fixed vertex is maximized around that distance.
The weight
\[
  J_i={d\choose i}
  \frac{\gamma^i(1-\gamma)^{d-i}}{\kappa_i}=
  (1-\gamma)^d\left\{\frac{1}{(1/\gamma-1)(n-1)}\right\}^i,
  \quad i\in \{1,\ldots,d\}
\]  
monotonically decreases in $i\in\{1,\ldots,d\}$ if
$\gamma\in(0,1-1/n)$. %, and increases if $\gamma\in(1-1/n,1)$.
When $\gamma=1-1/n$, this model reduces to the
distance-independent interaction model discussed in
Section~\ref{sect:model:ind}.

Since the analysis of this model is considerably
more involved than that of the two models discussed
above, we restrict our attention to the case $m>0$
and $\partial=\emptyset$.

We consider the random walk determined by 
\[
  w_i={d\choose i}\gamma^i(1-\gamma)^{d-i},
  \quad i\in\{0,1,\ldots,d\}.
\]
The eigenvalues \eqref{eigen} are 
$\lambda_i=\{1-n\gamma/(n-1)\}^i$, $i\in\{0,1,\ldots,d\}$.
This random walk is regarded as a long-range random
walk: the range $i\in\{0,1,\ldots,d\}$ is drawn from
the binomial distribution, and the next vertex is
chosen from the equally distant vertices with equal
probability.

The partition function \eqref{prop:part:z} becomes
\[
  Z_\mathcal{X}
  =\left(\frac{2\pi}{\beta}\right)^{n^d/2}\frac{1}{m}
  \prod_{i=1}^d
  \left\{1+m^2-\left(1-\frac{n\gamma}{(n-1)}\right)^i\right\}^{-\kappa_i/2},
\]
and the free energy converges to the distinct limits
as $d\to\infty$ or $n\to\infty$:
\begin{align*}
  &-\frac{1}{\beta}\lim_{d\to\infty}\frac{\log Z_{\mathcal{X}}}{n^d}
  =\frac{1}{2\beta}\left\{
  \log(1+m^2)-\log\left(\frac{2\pi}{\beta}\right)\right\};\\
  &-\frac{1}{\beta}\lim_{n\to\infty}\frac{\log Z_{\mathcal{X}}}{n^d}
  = \frac{1}{2\beta}\left[
    \log\{1+m^2-(1-\gamma)^d\}-\log\left(\frac{2\pi}{\beta}\right)
    \right],
\end{align*}
where the assumption of Corollary~\ref{cor:lim} is satisfied,
because
$\max_{i\in A_d(\epsilon)}|\lambda_i-\lambda_*|<2|1-n\gamma/(n-1)|^{(1-1/n)d}\to 0$ as $d\to \infty$, where $\lambda_*=\lim_{d\to\infty}\lambda_{\lfloor(1-1/n)d\rfloor}=0$.

The Green function \eqref{green_h} for this model is
\begin{equation}\label{cov_bin}
  G_{\alpha,\emptyset}(x,x')=
  \frac{1}{n^d}\sum_{i=0}^d
  \frac{K_i(\rho(x,x'))}{1-\alpha\left(1-\frac{n\gamma}{n-1}\right)^i}.
\end{equation}
The covariance decays exponentially with the Hamming distance
as $n\to\infty$, and we have \eqref{pro:limcov:1} with 
\begin{align*}
  C_{m,d}(\rho)&=
  \sum_{i=0}^\rho\binom{\rho}{i}
  \frac{(-1)^i}{1/\alpha-(1-\gamma)^{i+d-\rho}}
  =\alpha\sum_{i=0}^\rho\binom{\rho}{i}(-1)^i
  \sum_{k=0}^\infty\{\alpha(1-\gamma)^{i+d-\rho}\}^k
  \\
  &=\alpha\sum_{k=0}^\infty
  \left[\{\alpha(1-\gamma)^{d-\rho}\}^k
  \sum_{i=0}^\rho
  \binom{\rho}{i}(-1)^i(1-\gamma)^{ik}\right]\\
  &=\alpha\sum_{k=0}^\infty
  \{\alpha(1-\gamma)^{d-\rho}\}^k
  \{1-(1-\gamma)^k\}^\rho,
\end{align*}
where
$0 < C_{m,d}(\rho)\le\alpha\sum_{k=0}^\infty\alpha^k=\alpha/(1-\alpha)=1/m^2<\infty$.

We have the following estimate in the limit $d\to\infty$,
which is consistent with Theorem~\ref{thm:limcov}.
Here, we see that the convergence is fastest when
the weights are constant, i.e., when the model reduces
to the distance-independent interaction model.

\begin{proposition}\label{prop:binom}
  Set $m>0$ and $\partial=\emptyset$. For
  $(x,x')\in\mathcal{X}\times\mathcal{X}$,
  the Gaussian free field $(g_x)_{x\in\mathcal{X}}$
  defined via the Hamiltonian \eqref{ham_bin} satisfies
  \[
    {\rm Cov}(g_x,g_{x'})=
    \frac{\delta_{x,x'}}{\beta(1+m^2)}+
    O\{(1-\gamma)^d\},
    \quad d\to\infty,
  \]
  where $\gamma\in(0,1-1/n]$ and
  \[
    a_i(\gamma)=\left(1-\frac{1}{n}\right)
    \left(1-\frac{n\gamma}{n-1}\right)^i+\frac{1}{n} \in(0,1),
    \quad i\ge 0.
  \]
  In particular, the convergence of the variance
  to its limiting value is fastest when $\gamma=1-1/n$.
\end{proposition}

\begin{proof}
  Since $|\alpha(1-n\gamma/(n-1))|<1$, the expression
  \eqref{cov_bin} leads to
  \begin{align*}
    G_{\alpha,\emptyset}(x,x')&=\frac{1}{n^d}
    \sum_{i=0}^d\sum_{j=0}^\infty
    \left\{\alpha\Big(1-\frac{n\gamma}{n-1}\Big)^i\right\}^j
    K_i(\rho(x,x'))\\
    &=\frac{1}{n^d}
    \sum_{j=0}^\infty\alpha^j\sum_{i=0}^d
    \Big(1-\frac{n\gamma}{n-1}\Big)^{ji}
    K_i(\rho(x,x'))\\
    &=\frac{1}{(n-1)^{\rho(x,x')}}
    \sum_{j=0}^\infty\alpha^j
    \{a_j(\gamma)\}^{d-\rho(x,x')}\{1-a_j(\gamma)\}^{\rho(x,x')},
  \end{align*}
  where the third equality follows by \eqref{gen}. Then,
  we have
  \begin{align*}
    &G_{\alpha,\emptyset}(x,x)=1+
    \sum_{i=1}^\infty\alpha^i\{a_i(\gamma)\}^d, \\
  %\]
  %and  
  %\begin{equation}\label{prop:binom1}
    &G_{\alpha,\emptyset}(x,x')
    =\sum_{i=1}^\infty \alpha^i\{a_i(\gamma)\}^d
    \exp\left\{\rho(x,x')\log\frac{1-a_i(\gamma)}{a_i(\gamma)(n-1)}\right\},
    \quad x\neq x'.
  \end{align*}
  Since the sequence $a_i(\gamma)$, $i\ge 1$ is monotonically
  decreasing, we have
  \begin{align*}
    \alpha\{a_1(\gamma)\}^d=\alpha(1-\gamma)^d&<
    G_{\alpha,\emptyset}(x,x)-1
    =\sum_{i=1}^\infty\alpha^i\{a_i(\gamma)\}^d
    <\{a_1(\gamma)\}^d\sum_{i=1}^\infty\alpha^i\\
    &=(1-\gamma)^d\frac{\alpha}{1-\alpha},
   \end{align*}
  and we conclude that
  $|G_{\alpha,\emptyset}(x,x)-1|=\Theta\{(1-\gamma)^d\}$
  for large $d$.
  Thus, the convergence of the variance
  ${\rm Cov}(g_x,g_{x})=\alpha G_{\alpha,\emptyset}(x,x)/\beta$
  to its limiting value is fastest when $\gamma=1-1/n$. 
  For $x\neq x'$,
  %set $\rho(x,x')=sd$, $s\in(0,1]$.
  %For any small $\epsilon>0$, we can take large $J$ such that
  %$(s+n-1)\{1-n\gamma/(n-1)\}^J<\epsilon$. 
  %Bounding the summation \eqref{prop:binom1} from below by
  %the $J$-th term, we have
  %\[
  %  G_{\alpha,\emptyset}(x,x')>
  %  \alpha^Ja_J(\gamma)^d
  %  \left[\left\{\frac{1-a_J(\gamma)}{a_J(\gamma)(n-1)}\right\}^s
  %  \right]^d
  %  =\alpha^J\{a_J(\gamma)\}^d
  %  \left[\left\{\frac{1-a_J(\gamma)}{a_J(\gamma)(n-1)}\right\}^s
  %  \right]^d
  %\]
  \begin{align*}
    G_{\alpha,\emptyset}(x,x')
    <(1-\gamma)^d\sum_{i=1}^\infty \alpha^i
    \exp\left\{\rho(x,x')\log\frac{1-a_i(\gamma)}{a_i(\gamma)(n-1)}\right\}
    <(1-\gamma)^d\frac{\alpha}{1-\alpha},
  \end{align*}
  where the second inequality follows by 
  $0<\{1-a_i(\gamma)\}/\{a_i(\gamma)(n-1)\}<1$ for all $i\ge 1$.
  Therefore, we conclude that
  $|G_{\alpha,\emptyset}(x,x')|=O\{(1-\gamma)^d\}$ for large $d$,
  and we have
  ${\rm Cov}(g_x,g_{x'})=\alpha G_{\alpha,\emptyset}(x,x')/\beta=O\{(1-\gamma)^d\}$.
\end{proof}

Finally, let us consider a double-scaling limit, where
the interaction weights $J_i$, $i\in\{1,\ldots,d\}$
are scaled with the mass $m$ as $m,\gamma\to 0$ with
$\gamma/m^2\to 1-1/n$. In this limit, the correlation
between the scaled Gaussian random variables at two
vertices at the Hamming distance $\rho$ converges to
exactly the reciprocal of the number of vertices at
distance $\rho$ from a given vertex. In this sense, this
Gaussian free field completely encodes the geometry
of the Hamming graph.

\begin{theorem}\label{thm:binom:lim}
  Set $m^2>0$ and $\partial=\emptyset$. For
  $(x,x')\in\mathcal{X}\times\mathcal{X}$, the Gaussian free field
  $(g_x)_{x\in\mathcal{X}}$ defined via
  the Hamiltonian \eqref{ham_bin} satisfies
  \begin{align*}
    m^2n^d{\rm Cov}(g_x,g_{x'})\to &
    \frac{1}{\beta}\left[(d+1)\begin{pmatrix}d\\\rho\end{pmatrix}\left(1-\frac{1}{n}\right)^{\rho+1}\left(\frac{1}{n}\right)^{d-\rho}\right]^{-1}, \\
   & \quad m,\gamma\to 0,~\frac{\gamma}{m^2}\to 1-\frac{1}{n},
  \end{align*}
  where $\rho=\rho(x,x')$. In this limit, the correlation of the scaled field
  $(mg_x)_{x\in\mathcal{X}}$ satisfies
  $\lim_{m,\gamma\to 0}{\rm Corr}(mg_x,mg_{x'})=\kappa_{\rho}^{-1}$,
  where $\kappa_{\rho}$ is the number of vertices at the Hamming
  distance $\rho$ from a given vertex.
\end{theorem}  
  
\begin{proof}
Setting $n\gamma/(n-1)=m^2$, the eigenvalues become $\lambda_i=(1-m^2)^i$,
and we have
\begin{align*}
  &\alpha G_{\alpha,\emptyset}(x,x')=
  \frac{1}{(n-1)^\rho}
  \sum_{j=0}^\infty\left(\frac{1}{1+m^2}\right)^{j+1}
  \left[1-\left(1-\frac{1}{n}\right)\left\{1-(1-m^2)^j\right\}\right]^{d-\rho}\\
  &\times
  \left[\left(1-\frac{1}{n}\right)\left\{1-(1-m^2)^j\right\}\right]^{\rho}\\
  &=
  \frac{1}{(n-1)^\rho}
  \sum_{j=0}^\infty\left(\frac{1}{1+m^2}\right)^{j+1}\sum_{k=0}^{d-\rho}
  \begin{pmatrix}d-\rho\\k\end{pmatrix}(-1)^k
  \left[\left(1-\frac{1}{n}\right)\left\{1-(1-m^2)^j\right\}\right]^{k+\rho}\\
  &=
  \frac{1}{n^\rho}
  \sum_{j=0}^\infty\left(\frac{1}{1+m^2}\right)^{j+1}\sum_{k=0}^{d-\rho}
  \begin{pmatrix}d-\rho\\k\end{pmatrix}
    \left(\frac{1}{n}-1\right)^k
  \sum_{l=0}^{k+\rho}\begin{pmatrix}k+\rho\\l\end{pmatrix}
    \{-(1-m^2)^j\}^{l},\\
  &=
  \frac{1}{n^\rho}
  \sum_{k=0}^{d-\rho}\sum_{l=0}^{k+\rho}(-1)^l
  \begin{pmatrix}k+\rho\\l\end{pmatrix}
  \begin{pmatrix}d-\rho\\k\end{pmatrix}
    \left(\frac{1}{n}-1\right)^k
    \sum_{j=0}^\infty
    \frac{(1-m^2)^{jl}}{(1+m^2)^{j+1}}\\
  &=
  \frac{1}{n^\rho}
  \sum_{k=0}^{d-\rho}\sum_{l=0}^{k+\rho}(-1)^l
  \begin{pmatrix}k+\rho\\l\end{pmatrix}
    \begin{pmatrix}d-\rho\\k\end{pmatrix}
    \left(\frac{1}{n}-1\right)^k
    \frac{1}{1+m^2-(1-m^2)^l}.
\end{align*}
For small $m^2$, we have
\begin{align*}
  &\alpha G_{\alpha,\emptyset}(x,x')=
  \frac{1}{m^2n^\rho}
  \sum_{k=0}^{d-\rho}\sum_{l=0}^{k+\rho}
  \frac{(-1)^l}{l+1}
  \begin{pmatrix}k+\rho\\l\end{pmatrix}
    \begin{pmatrix}d-\rho\\k\end{pmatrix}
    \left(\frac{1}{n}-1\right)^k
    \{1+O(m^2)\}\\
  &=
  \frac{1}{m^2n^\rho}
  \sum_{k=0}^{d-\rho}\sum_{l=1}^{k+\rho+1}
  \begin{pmatrix}k+\rho+1\\l\end{pmatrix}(-1)^{l-1}
    \begin{pmatrix}d-\rho\\k\end{pmatrix}
    \left(\frac{1}{n}-1\right)^k
    \frac{1}{k+\rho+1}+O(1)\\
  &=
  \frac{1}{m^2n^\rho}
  \sum_{k=0}^{d-\rho}
    \begin{pmatrix}d-\rho\\k\end{pmatrix}
    \left(\frac{1}{n}-1\right)^k
    \frac{1}{k+\rho+1}+O(1)\\
 &=
  \frac{1}{m^2n^\rho}\int_0^1u^\rho
  \sum_{k=0}^{d-\rho}
    \begin{pmatrix}d-\rho\\k\end{pmatrix}
   \left\{\left(\frac{1}{n}-1\right)u\right\}^kdu+O(1)\\
&=
  \frac{1}{m^2n^\rho}\int_0^1u^\rho
  \left\{1-\left(1-\frac{1}{n}\right)u\right\}^{d-\rho}du+O(1)\\
&=\left[m^2n^d(d+1)\begin{pmatrix}d\\\rho\end{pmatrix}\left(1-\frac{1}{n}\right)^{\rho+1}\left(\frac{1}{n}\right)^{d-\rho}\right]^{-1}+O(1).
\end{align*}
Therefore, the covariance is
\begin{align*}
m^2n^d{\rm Cov}(g_x,g_{x'})&=
m^2n^d\frac{\alpha}{\beta}G_{\alpha,\emptyset}(x,x')\\
&=
\frac{1}{\beta}\left[(d+1)\begin{pmatrix}d\\\rho\end{pmatrix}\left(1-\frac{1}{n}\right)^{\rho+1}\left(\frac{1}{n}\right)^{d-\rho}\right]^{-1}+O(1).
\end{align*}
The correlation between the scaled Gaussian random
variables at the Hamming distance $\rho$
converges to the limit
\[
\lim_{m,\gamma\to 0}{\rm Corr}\left(mg_x,mg_{x'}\right)=
\frac{1}{n^d}
  \left[\begin{pmatrix}d\\\rho\end{pmatrix}\left(1-\frac{1}{n}\right)^{\rho}\left(\frac{1}{n}\right)^{d-\rho}\right]^{-1}
=\frac{1}{\kappa_\rho}.
\]
\end{proof}

Theorem~\ref{thm:binom:lim} showed that the correlation
of the scaled field $(mg_x)_{x\in\mathcal{X}}$, viewed
as a function of the Hamming distance, is convex and
achieves its minimum at a point strictly in the interior
of $\{0,1,\ldots,d\}$. This is because $\kappa_\rho/n^d$
is the probability mass function of a binomial distribution
whose mode is around $d(1-1/n)$. This does not contradict
the assertion of Corollary~\ref{cor:cov}, as that corollary
assumes the interaction weights and the mass $m$ are
independent. Counter-intuitively, the correlation
with the most distant vertex can be stronger than with
adjacent one; for $n=2$, for instance, the correlation
between diagonally opposite vertices exceeds that
between adjacent vertices.

\section*{Acknowledgements}
The author was supported in part by JSPS KAKENHI Grants
24K06876. The author is grateful to Robert Griffiths for
fruitful discussions. The author also thanks the anonymous
referees for their careful reading of the manuscript and
helpful comments.

%\section*{Appendix}

% To print the credit authorship contribution details
%\printcredits

%% Loading bibliography style file
%\bibliographystyle{model1-num-names}
%\bibliographystyle{cas-model2-names}

% Loading bibliography database
%\bibliography{cas-refs}

% Biography
%\bio{}
% Here goes the biography details.
%\endbio

%\bio{pic1}
% Here goes the biography details.
%\endbio

\begin{flushleft}

Shuhei Mano\\
The Institute of Statistical Mathematics, Tokyo 190-8562, Japan\\
E-mail: smano@ism.ac.jp

\end{flushleft}

\end{document}